%% file: AsDist.tex
\documentclass[11pt,a4paper]{article}
\usepackage{amsmath,amsfonts,latexsym,amssymb,chicago,annals3,graphicx,here}
\advance \topmargin by -\headheight \advance \topmargin by
-\headsep

\def\P{{\mathbb P}}

\def\H{{\mathbb H}}
\def\t{\tau}

\input commands3.tex

\begin{document}
\title{Testing monotonicity of a hazard: asymptotic distribution theory}
\author{Piet Groeneboom and Geurt Jongbloed}
\date{\today}
\affiliation{Delft University of Technology}
\AMSsubject{Primary: 62E20, secondary: 62N03.}
\maketitle

\begin{abstract} Two new test statistics are introduced to test the null hypotheses that the sampling distribution has an increasing hazard rate on a specified interval $[0,a]$. These statistics are empirical $L_1$-type distances between the isotonic estimates, which use the monotonicity constraint, and either the empirical distribution function or the empirical cumulative hazard. They measure the excursions of the empirical estimates with respect to the isotonic estimates, due to local non-monotonicity. Asymptotic normality of the test statistics, if the hazard is strictly increasing on $[0,a]$, is established under mild conditions. This is done by first approximating the global empirical distance by an distance with respect to the underlying distribution function. The resulting integral is treated as sum of increasingly many local integrals to which a CLT can be applied. The behavior of the local integrals is determined by a canonical process: the difference between the stochastic process $x\mapsto W(x)+x^2$ where $W$ is standard two-sided Brownian Motion, and its greatest convex minorant.
\end{abstract}

\section{Introduction}
\label{sec:intro}
One way of characterizing a distribution of an absolutely continuous random variable $X$ that is particularly useful in reliability theory and survival analysis, is by its hazard rate $h_0$. Suppose $X$ models the failure time of a certain device. The interpretation of the hazard rate is that for small $\epsilon>0$, $\epsilon h_0(x)$ reflects the probability of failure of the device in the time interval $(x,x+\epsilon]$ given the device was still unimpaired at time $x$ (assuming $h_0$ is continuous at $x$). Put differently, $h_0(x)$ represents the level of instantaneous risk of failure of the device at  time $x$, given it still works at time $x$. A high value reflects high risk, a low value low risk. Lifetimes of devices that are subject to aging can be described by distributions with increasing hazard rate. Locally decreasing hazard rates can be used to model life times of devices that become more reliable with age during a certain period of time.

It is especially this clear interpretation of these qualitative properties of a hazard rate that makes this function a natural characteristic of a survival distribution.  The problem of estimating a hazard rate nonparametrically under qualitative (or shape) restrictions gained attention in the sixties of the previous century (see \mycite{GrJo11est} and the references therein). Also the problem of testing the null hypothesis of constant hazard (exponentiality) against monotonicity of the hazard was studied intensively, see e.g.\ \mycite{prospyke:67}. Only quite recently another testing problem, with a ``shape constraint'' rather than parametric null hypothesis was studied. See also the discussion in companion paper \mycite{GrJo10}.

In this paper, we consider the asymptotic distribution theory for two  integral-type test statistics for the hypothesis that a hazard rate $h_0$ is  monotone on an interval $[0,a]$, for some known $a>0$. We restrict ourselves to the increasing case; the case of locally decreasing hazard can be considered analogously.
%
%

Based on an i.i.d.\ sample $X_1,\ldots,X_n$ from the distribution associated with $H_0$, the most natural nonparametric estimator for $H_0$ without assuming anything on $H_0$, is  the empirical cumulative hazard function given by
$$
\H_n(x)=\left\{
\begin{array}{lll}
-\log\left\{1-\F_n(x)\right\},\,&x\in\left[0,X_{(n)}\right),\\
\infty,\,&x\ge X_{(n)}
\end{array}
\right.
$$
where $\F_n$ denotes the empirical distribution function based on $X_1,X_2,\ldots,X_n$. Under the assumption that $H_0$ is convex on $[0,a]$,  the cumulative hazard can be estimated by the greatest convex minorant $\hat H_n$ of the empirical cumulative hazard function $\H_n$ on the interval $[0,a]$. Using these two estimators, the following test statistic emerges:
\begin{equation}
\label{eq:defTn}
T_n=\int_{[0,a]}\bigl\{\H_n(x-)-\hat H_n(x)\bigr\}\,d\F_n(x).
\end{equation}
Note that this is the empirical $L_1$-distance  between the two mentioned estimators for the cumulative hazard function w.r.t.\ the empirical measure $d\F_n$, and that $T_n\ge0$ since $\hat H_n$ is a minorant of $\H_n$. If $H_0$ is concave on $[0,a]$, both estimators for $H_0$ will be close to $H_0$ and $T_n$ will tend to be small (converge to zero a.s.\  for $n\rightarrow\infty$). On the contrary, if $h_0$ has a region in $[0,a]$ where it is not increasing, $\H_n$ will capture this ``non-convexity" of $H_0$ and converge to $H_0$ on this region whereas $\hat H_n$ will converge to the convex minorant of $H_0$ on $[0,a]$. Note that $T_n=0$ if and only if $\hat{H}_n$ coincides with the linear interpolation of the points $(x_{(i)},\H_n(x_{(i)}-))$ on the range of the data falling in $[0,a]$. One could say that $T_n=0$ if $\H_n$ is `as convex as it can be on $[0,a]$', being an increasing right continuous step function. This is the reason for taking $\H_n(x-)$ instead of $\H_n(x)$ in (\ref{eq:defTn}). A similar reasoning can be held for another test statistic,
\begin{equation}
\label{eq:defUn}
U_n=\int_{[0,a)}\bigl\{\F_n(x-)-\hat F_n(x)\bigr\}\,d\F_n(x), \mbox{ where } \hat{F}_n(x)=1-\exp(-\hat{H}_n(x)).
\end{equation}
An advantage of this definition is that $U_n$ is less sensitive to possible problems that can occur with large values of $\H_n$.

The main result of this paper concerns the asymptotic distribution of $T_n$ and $U_n$: under certain assumptions
\begin{equation}
\label{eq:limres}
n^{5/6}\left\{T_n-E\tilde T_n\right\}\stackrel{{\cal D}}\longrightarrow N\left(0,\sigma_{H_0}^2\right)
\mbox{ and }n^{5/6}\left\{U_n-EU_n\right\}\stackrel{{\cal D}}\longrightarrow N\left(0,\sigma_{F_0}^2\right),
\end{equation}
where $\tilde T_n$ is a modified version of $T_n$, see Theorems \ref{th:hazard_theorem} and  \ref{th:df_theorem}. Here $\sigma_{H_0}^2$ and $\sigma_{F_0}^2$ are constants depending on $f_0$. Results of a similar flavor were established in, e.g., \mycite{vladrik:08} for the difference between the empirical distribution function and its concave majorant.

The basic idea of the proof  is to approximate the integral in the test statistic as sum of increasingly many local integrals, using the crucial localization Lemma \ref{localization_lemma}, and to apply a Central Limit Theorem to the components that arise in this way. The behavior of the local integrals is determined by a canonical process, the difference between a Brownian motion with parabolic drift and its convex minorant. Relevant properties of this process are derived in section \ref{sec:asproc}. In section \ref{sec:embedding}, a statistic related to $T_n$ (where the integral is taken w.r.t.\ $F_0$ rather than $\F_n$) is closely approximated by an integral involving the independent increments of Brownian motion. Moreover, the resulting integral is represented as a sum of local integrals using a ``big blocks separated by small blocks'' construction as introduced in \mycite{rosenblatt:56}. The local integrals over the big blocks reduce to the processes considered in section \ref{sec:asproc}. Finally, because the local integrals are based on the independent increments of a Brownian motion process, a CLT can be applied to obtain the first result in the spirit of (\ref{eq:limres}). In section  \ref{sec:further}, the main results of the paper are established by showing that the differences between the integrals w.r.t.\ $d\F_n$ and $dF_0$ are sufficiently small.

\section{Asymptotic local problem}
\label{sec:asproc}
\begin{figure}[!ht]
\begin{center}
\includegraphics[scale=0.5]{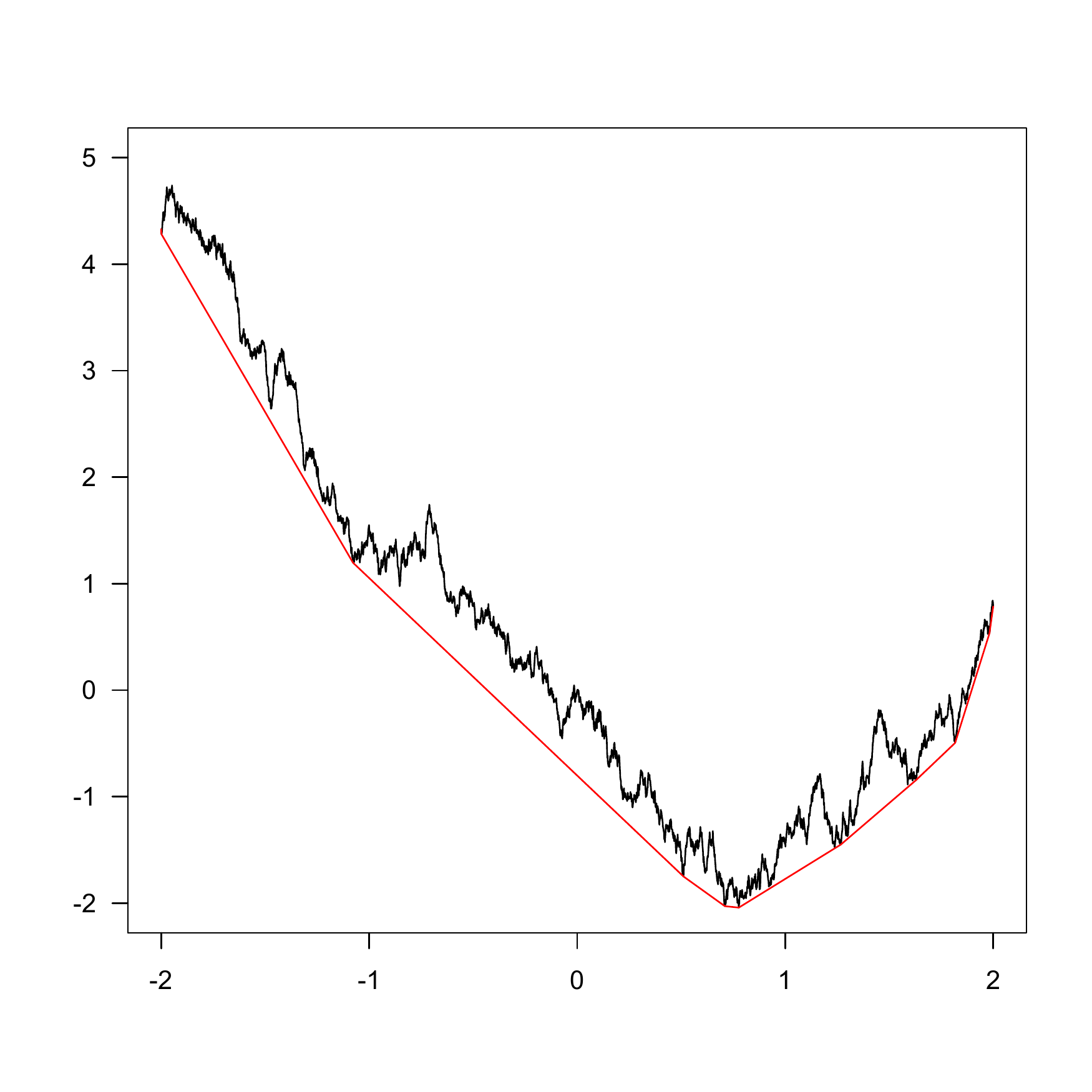}
\end{center}
\caption{The greatest convex minorant of $W(x)+x^2$, restricted to $[-2,2]$.}
\label{fig:convex_min}
\end{figure}

Consider the process
\begin{equation}
\label{eq:defV}
x\mapsto V(x)=W(x)+x^2,\,x\in\R
\end{equation}
with $W$ standard two-sided Brownian motion on $\R$.  Then, for $c>0$, define the functional  $Q_c$ by
\begin{equation}
\label{Q_c}
Q_c=\int_0^c\left\{V(x)-C(x)\right\}\,dx,
\end{equation}
where  $C$ is the greatest convex minorant of $V$ on $\R$. For a picture of the process $V$ and its greatest convex minorant, restricted to the interval $[-2,2]$, see Figure \ref{fig:convex_min}. We have the following result.

\begin{theorem}
\label{th:BM-limit}
$$
c^{-1/2}\left\{Q_c-c E|C(0)|\right\}\stackrel{{\cal D}}\longrightarrow N(0,\s^2),\,c\to\infty,
$$
where $C(0)$ is the value of the greatest convex minorant $C$ of the process $V$ at zero, and
$$
\s^2=2\int_0^{\infty}\mbox{\rm covar}(-C(0),V(x)-C(x))\,dx.
$$
All moments of $c^{-1/2}\left\{Q_c-c E|C(0)|\right\}$ exist and (in particular) the fourth moment is uniformly bounded in $c$ and converges to the fourth moment of the normal $N(0,\s^2)$ distribution, as $c\to\infty$.
\end{theorem}

In the proof we will use the following lemma, which is proved in the appendix.

\begin{lemma}
\label{lem:Vbelow0}
For the process $V$ defined in (\ref{eq:defV}), there exist positive constants $c$ and $c^\prime$ such that for all $u\ge0$
$$
P\left(\min_{x\not\in [-u,u]}V(x)\le0\right)\le c e^{-c^\prime u^3}.
$$
\end{lemma}

\noindent
{\bf Proof of Theorem \ref{th:BM-limit}.}
It follows from the results in \mycite{piet:89} that the process
\begin{equation}
\label{stationarity1}
V(x)-C(x),\,x\in\R,
\end{equation}
is stationary. In fact, the process touches zero at changes of slope of $C$ and behaves between these touches of zero as an excursion of a Brownian motion path above a parabola of the form
$$
\f(x)=s-(x-a)^2,\,x\in\R,
$$
where $\f$ is a parabola touching two local minima of Brownian motion, and where the (random) values $a$ and $s$ depend on the Brownian motion path. Defining
$$
D_k=\int_{k}^{k+1}\left\{V(x)-C(x)\right\}\,dx,\,k\in{\mathbb Z},
$$
we get a stationary sequence of random variables, and the stationarity of the process (\ref{stationarity1})
yields:
$$
E D_k=\int_{k}^{k+1}E\left\{V(x)-C(x)\right\}\,dx=E|C(0)|.
$$

Moreover, all moments of $D_k$ exist. This follows from the fact that
$$
\max_{x\in[0,1]}\{V(x)-C(x)\}
$$
has a distribution with tails which die out faster than exponentially. To see this, note that, $\forall u\ge0$,
\begin{align}
&\P\left\{\max_{x\in[0,1]}\{V(x)-C(x)\}\ge M\right\}
\le \P\left\{\max_{x\in[0,1]}V(x)\ge \tfrac12M\right\}+\P\left\{\min_{x\in\R}C(x)\le -\tfrac12M\right\}\nonumber\\
&\le \P\left\{\max_{x\in[0,1]}W(x)\ge \tfrac12M-1\right\}+\P\left\{\min_{x\in\R}V(x)\le -\tfrac12M\right\}\nonumber\\
&\le\sqrt{\frac2{\pi}}\int_{\tfrac12M-1}^{\infty}e^{-\tfrac12x^2}\,dx+\P\left\{\min_{x\in[-u,u]}W(x)\le -\tfrac12M\right\}
+\P\left\{\min_{x\notin[-u,u]}V(x)\le 0\right\}.\label{eq:UBp}
\end{align}

The first term on the right hand side is bounded by $c\exp\{-c^{\prime}M^2/4\}$ for some $c,c^{\prime}>0$. By Lemma \ref{lem:Vbelow0} we have for the third term:
$$
P\left(\min_{x\not\in [-u,u]}V(x)\le0\right)\le c e^{-c^\prime u^3}
$$
for constants $c,c'>0$. For the second term in (\ref{eq:UBp}),
we get by Brownian scaling,
\begin{align*}
&\P\left\{\min_{x\in[-u,u]}W(x)\le -M/2\right\}=\P\left\{\min_{x\in[-1,1]}W(ux)\le -M/2\right\}\\
&=\P\left\{\min_{x\in[-1,1]}u^{-1/2}W(ux)\le -u^{-1/2}M/2\right\}
=\P\left\{\max_{x\in[-1,1]}W(x)\ge u^{-1/2}M/2\right\}\\
&\le2\sqrt{\frac2{\pi}}\int_{u^{-1/2}M/2}^{\infty}e^{-\tfrac12x^2}\,dx\le2\sqrt{\frac2{\pi}}\frac{2\sqrt{u}}{M}\exp\left\{-\tfrac18M^2/u\right\}.
\end{align*}
Hence, taking $u=\sqrt{M}$ in the second and third term  in (\ref{eq:UBp}) and observing that the first term is of lower order, we obtain
\begin{equation}
\label{moment_condition}
\P\{\max_{x\in[0,1]}\{V(x)-C(x)\}\ge M\}\le  c_1e^{-c_2M^{3/2}}
\end{equation}
for constants $c_1,c_2>0$.

Now let $\t(a)$ be defined by:
$$
\t(a)=\mbox{argmin}_{x\in\R}\left\{W(x)+(x-a)^2\right\}.
$$
The (stationary) process $a\mapsto \t(a)-a$ is studied in \mycite{piet:89} and it follows from the results, given there, that there exist positive constants $c_1$ and $c_2$ such that
$$
|\P(A\cap B)-\P(A)\P(B)|\le c_1e^{-c_2m^3},
$$
for events $A$ and $B$ such that
$$
A\in\s\left\{\t(a):a\le0\right\},\qquad B\in\s\left\{\t(a):a\ge m\right\}.
$$
This implies that there also exist positive constants $c_1$ and $c_2$ such that
\begin{equation}
\label{mixing_condition}
|\P(A\cap B)-\P(A)\P(B)|\le c_1e^{-c_2m^3},
\end{equation}
for events $A$ and $B$ such that
$$
A\in\s\left\{V(x)-C(x):x\le0\right\},\qquad B\in\s\left\{V(x)-C(x):x\ge m\right\}.
$$
So we can apply Theorem 18.5.3 in \mycite{ibralin:71}, p.\ 347, yielding that
$$
c^{-1/2}\left\{Q_c-c E|C(0)|\right\}\stackrel{{\cal D}}\longrightarrow N(0,\s^2),
$$
where
$$
\s^2={\rm var}(D_0)+2\sum_{k=1}^{\infty}{\rm covar}(D_0,D_k).
$$
Using the stationarity of the process (\ref{stationarity1}) again, we obtain
$$
\s^2=2\int_0^{\infty}\mbox{\rm covar}(-C(0),V(x)-C(x))\,dx.
$$
The last statement of the theorem follows from (\ref{moment_condition}) and (\ref{mixing_condition}).
\eop

\vspace{0.3cm}
We will also need the following extension of Theorem \ref{th:BM-limit}.

\begin{theorem}
\label{th:BM-limit2}
Let $C_c$ be the greatest convex minorant on $[0,c]$ of the process
$$
V(x),\,x\in[0,c].
$$
Note that $C_c$ is not the restriction of $C$ to $[0,c]$, since $C$ is globally defined on $\R$, and $C_c$ is the greatest convex minorant of the process $V$ on $[0,c]$, and only defined on $[0,c]$.
\begin{enumerate}
\item[(i)] Let, for $c>4$, the interval $I_c$ be defined by
$$
I_c=\left[\sqrt{c},c-\sqrt{c}\right].
$$
Then:
\begin{equation}
\label{restricted_CM}
c^{-1/2}\left\{\int_{I_c}\left\{V(x)-C_c(x)\right\}\,dx-E\int_{I_c}\left\{V(x)-C_c(x)\right\}\,dx\right\}\stackrel{{\cal D}}\longrightarrow N(0,\s^2),\,c\to\infty,
\end{equation}
where $\s^2$ is defined as in Theorem \ref{th:BM-limit}.
\item[(ii)] Relation (\ref{restricted_CM}) also holds if the interval $I_c$ is given by:
$$
I_c=\left[0,c-\sqrt{c}\right]\mbox{ or }I_c=\left[\sqrt{c},c\right].
$$
\item[(iii)] For any choice of $I_c$ in (i) or (ii), the fourth moment of
$$
c^{-1/2}\left\{\int_{I_c}\left\{V(x)-C_c(x)\right\}\,dx-E\int_{I_c}\left\{V(x)-C_c(x)\right\}\,dx\right\}
$$
is uniformly bounded in $c$, and converges to the fourth moment of a normal $N(0,\s^2)$ distribution, as $c\to\infty$.
\end{enumerate}
\end{theorem}

\noindent
{\bf Proof.} (i). The probability that $C_c$ is different from $C$ on the interval $I_c$ is less than or equal to
$$
k_1\exp\left\{-k_2c^{3/2}\right\},
$$
for constants $k_1,k_2>0$. The proof of this is analogous to the proof of Lemma \ref{localization_lemma} in the next section. Hence, if $K_c$ denotes the event that $C_c\not\equiv C$ on $I_c$, we get:
$$
E\int_{I_c} \left|C_c(x)-C(x)\right|\,dx
\le \left\{\int_{I_c} E\left\{V(x)-C(x)\right\}^2\,dx\right\}^{1/2}\P(K_c)^{1/2}
=O\left(c^{1/2}e^{-kc^{3/2}}\right),\,c\to\infty,
$$
for some $k>0$. Hence:
\begin{align*}
&c^{-1/2}\left\{\int_{I_c}\left\{V(x)-C_c(x)\right\}\,dx-E\int_{I_c}\left\{V(x)-C_c(x)\right\}\,dx\right\}\\
&=c^{-1/2}\left\{\int_{I_c}\left\{V(x)-C(x)\right\}\,dx-E\int_{I_c}\left\{V(x)-C(x)\right\}\,dx\right\}+O_p\left(ce^{-kc^{3/2}}\right),
\end{align*}
and the statement now follows.\\
(ii). We can repeat the argument on the interval $[0,\sqrt{c}]$, and apply the argument used in (i) on the subinterval $I'_c=[c^{1/4},\sqrt{c}-c^{1/4}]$ (but leaving $C_c$ as it was defined in (i)). This yields:
$$
c^{-1/4}\left\{\int_{I_c'}\left\{V(x)-C_c(x)\right\}\,dx-E\int_{I_c'}\left\{V(x)-C_c(x)\right\}\,dx\right\}\stackrel{{\cal D}}\longrightarrow N(0,\s^2),\,c\to\infty,
$$
implying:
$$
c^{-1/2}\left\{\int_{I_c'}\left\{V(x)-C_c(x)\right\}\,dx-E\int_{I_c'}\left\{V(x)-C_c(x)\right\}\,dx\right\}\stackrel{p}\longrightarrow 0,\,c\to\infty,
$$
Moreover,
$$
c^{-1/2}\int_{[0,c^{1/4}]}E\left|V(x)-C_c(x)\right|\,dx=O\left(c^{-1/4}\right),\,c\to\infty.
$$
The statement now follows for the first choice of the interval $I_c$ in (ii). For the second choice of $I_c$ the argument is similar.\\
(iii). Let $I_c$ be as in (i). Then:
\begin{align*}
&c^{-2}E\left\{\int_{I_c}\left\{V(x)-C_c(x)\right\}\,dx-E\int_{I_c}\left\{V(x)-C_c(x)\right\}\,dx\right\}^4\\
&=c^{-2}E\left\{\int_{I_c}\left\{V(x)-C(x)\right\}\,dx-E\int_{I_c}\left\{V(x)-C(x)\right\}\,dx\right\}^4+O\left(e^{-kc^{3/2}}\right),
\end{align*}
for some $k>0$, and the statement now follows from Theorem \ref{th:BM-limit}, (\ref{moment_condition}) and (\ref{mixing_condition}) and the fact that
$$
\left(\frac{c}{c-2\sqrt{c}}\right)^2\to1,\,c\to\infty.
$$
If, for example, $I_c=[0,c-\sqrt{c}]$, we write
\begin{align*}
\int_{I_c}\left\{V(x)-C_c(x)\right\}\,dx-E\int_{I_c}\left\{V(x)-C_c(x)\right\}\,dx
=A_c+B_c,
\end{align*}
where
$$
A_c=\int_{[0,\sqrt{c}]}\left\{V(x)-C_c(x)\right\}\,dx-E\int_{[0,\sqrt{c}]}\left\{V(x)-C_c(x)\right\}\,dx
$$
and
$$
B_c=\int_{[\sqrt{c},c-\sqrt{c}]}\left\{V(x)-C_c(x)\right\}\,dx-E\int_{[\sqrt{c},c-\sqrt{c}]}\left\{V(x)-C_c(x)\right\}\,dx.
$$
Hence we get:
\begin{align*}
&c^{-2}E\left\{\int_{I_c}\left\{V(x)-C_c(x)\right\}\,dx-E\int_{I_c}\left\{V(x)-C_c(x)\right\}\,dx\right\}^4\\
&=c^{-2}EB_c^4+c^{-2}\left\{4EB_c^3A_c+6EB_c^2A_c^2+4EB_cA_c^3+EA_c^4\right\}.
\end{align*}
We have:
$$
c^{-2}EA_c^4=\left(\frac{\sqrt{c}}{c}\right)^2c^{-1}EA_c^4=O\left(c^{-1}\right),
$$
and similarly, using the Cauchy-Schwarz inequality,
$$
c^{-2}EB_cA_c^3=Ec^{-3/2}A_c^3 c^{-1/2}B_c\le\sqrt{Ec^{-3}A_c^6}\sqrt{Ec^{-1}B_c^2}=O\left(c^{-3/4}\right).
$$
Continuing in this way, we find that the only non-vanishing term is $c^{-2}EB_c^4$.
The statement now follows from what we proved for $I_c=[\sqrt{c},c-\sqrt{c}]$.\eop

\vspace{0.3cm}
We finally also need the following extension of Theorem \ref{th:BM-limit2}.

\begin{theorem}
\label{th:BM-limit3}
Let $F_c$, $G_c$ and $H_c$ be twice differentiable increasing functions on $[0,c]$, with continuous derivatives $f_c$, $g_c$ and $h_c$, respectively, satisfying
$$
F_c(x)= f_c(0)x(1+o(1)),\qquad G_c(x)= g_c(0)x(1+o(1)),\qquad H_c(x)=\tfrac12 h_c'(0)x^2(1+o(1)),\,c\to\infty,
$$
where the $o(1)$ term is uniform in $x$. We assume that $f_c(0)$, $g_c(0)$, $h_c(0)$ and $h_c'(0)$ are positive and stay away from zero and $\infty$, as $c\to\infty$, where $h_c'(0)$ denotes the right derivative of $h_c$ at zero.
Let $C_c$ be the greatest convex minorant on $[0,c]$ of the process
$$
V_c(x)=H_c(x)+W(G_c(x)),\,x\in[0,c].
$$
Moreover, let $S_c$ be defined by
$$
S_c(x)=V_c(x)-C_c(x),\,x\in[0,c].
$$
Then:
\begin{enumerate}
\item[(i)] Let, for $c>4$, the interval $I_c$ be defined by
$$
I_c=\left[\sqrt{c},c-\sqrt{c}\right].
$$
Then:
\begin{equation}
\label{var_approximation}
c^{-1}E\int_{I_c}S_c(x)\,dF_c(x)\sim\frac{g_c(0)^{2/3}f_c(0)}{\left(\tfrac12 h_c'(0)\right)^{1/3}}E|C(0)|,\,
\mbox{\rm var}\left(c^{-1/2}\int_{I_c}S_c(x)\,dF_c(x)\right)\sim \s_c^2,\,c\to\infty,
\end{equation}
where
\begin{equation}
\label{sigma_c}
\s_c^2=\frac{g_c(0)^{5/3}f_c(0)^2}{\left(\tfrac12 h_c'(0)\right)^{4/3}}
\sigma^2,
\end{equation}
and  $C$ and $\sigma^2$ are defined as in Theorem \ref{th:BM-limit}. Moreover, the fourth moment of
$$
c^{-1/2}\int_{I_c}\left\{S_c(x)-ES_c(x)\right\}\,dF_c(x)
$$
is uniformly bounded, and satisfies:
\begin{equation}
\label{4th_moment}
E\left(c^{-1/2}\int_{I_c}\left\{S_c(x)-ES_c(x)\right\}\,dF_c(x)\right)^4\sim M_c^{(4)},\,c\to\infty,
\end{equation}
where $M_c^{(4)}$ denotes the fourth moment of a normal $N(0,\s_c^2)$ distribution.
\item[(ii)] Relations (\ref{var_approximation}) and (\ref{4th_moment}) also hold if the interval $I_c$ is given by:
$$
I_c=\left[0,c-\sqrt{c}\right],\,I_c=\left[\sqrt{c},c\right]\mbox{ or }I_c=[0,c].
$$
\end{enumerate}
\end{theorem}

\noindent
{\bf Proof.} Since the proof proceeds along the lines of the proofs of Theorems \ref{th:BM-limit} and \ref{th:BM-limit2}, we only pay attention to the new type of scaling which is present in the process
$$
x\mapsto \tfrac12 h_c^\prime(0)x^2+W(g_c(0)x),\,x\in[0,c],
$$
which replaces the process
$$
x\mapsto V(x)=x^2+W(x),\,x\in[0,c].
$$

Let $a,b>0$. By Brownian scaling, the process
\begin{equation}
\label{eq:scale1}
x\mapsto ax^2+W(bx),\,x\in[0,c],
\end{equation}
has the same distribution as the process
\begin{equation}
\label{eq:scale2}
x\mapsto a^{-1/3}b^{2/3}\left\{\left(a^{2/3}b^{-1/3}x\right)^2+W(a^{2/3}b^{-1/3}x)\right\},\, x\in[0,c].
\end{equation}
Hence, if $C_{a,b}$ is the greatest convex minorant of the process given in (\ref{eq:scale1})
and $\tilde{C}_{a,b}$ of the process given in (\ref{eq:scale2})
we get:
\begin{align*}
&\int_0^c\left\{ax^2+W(bx)-C_{a,b}(x)\right\}f_c(0)\,dx\\
&  \=in^{\cal D} a^{-1/3}b^{2/3}f_c(0)\int_0^c\left\{\left(a^{2/3}b^{-1/3}x\right)^2+W(a^{2/3}b^{-1/3}x)
-a^{1/3}b^{-2/3}\tilde{C}_{a,b}(x)\right\}\,dx\\
&=\frac{bf_c(0)}{a}\int_0^{a^{2/3}b^{-1/3}c}\left\{u^2+W(u)
-a^{1/3}b^{-2/3}\tilde{C}_{a,b}\left(a^{-2/3}b^{1/3}u\right)\right\}\,du\\
&=\frac{bf_c(0)}{a}\int_0^{a^{2/3}b^{-1/3}c}\left\{u^2+W(u)
-C_c(u)\right\}\,du,
\end{align*}
where $C_c$ is the greatest convex minorant of the process
$$
u\mapsto u^2+W(u),\,u\in \left[0,a^{2/3}b^{-1/3}c\right].
$$

Thus, for $c\rightarrow\infty$
\begin{align*}
c^{-1}E\int_0^c\left\{ax^2+W(bx)-C_{a,b}(x)\right\}f_c(0)\,dx
\sim \frac{b^{2/3}f_c(0)}{a^{1/3}}E |C(0)|.
\end{align*}
Using that $a=\tfrac12 h_c'(0)$, $b=g_c(0)$, (\ref{var_approximation}) follows. Moreover,
\begin{align*}
&\mbox{var}\left(c^{-1/2}\int_0^c\left\{ax^2+W(bx)-C_{a,b}(x)\right\} f_c(0)\,du\right)\\
&=\frac{b^2f_c(0)^2}{a^2c}\mbox{var}\left(\int_0^{a^{2/3}b^{-1/3}c}\left\{u^2+W(u)
-C_c(u)\right\}\,du\right)\\
&=\frac{b^{5/3}f_c(0)^2}{a^{4/3}}\mbox{var}\left(\frac1{\sqrt{b^{-1/3}a^{2/3}c}}\int_0^{a^{2/3}b^{-1/3}c}\left\{u^2+W(u)
-C_c(u)\right\}\,du\right)\\
&\sim\frac{b^{5/3}f_c(0)^2}{a^{4/3}}\sigma^2\,c\to\infty.
\end{align*}
Taking $a=\tfrac12 h_c'(0)$, $b=g_c(0)$ now yields (\ref{sigma_c}).
\eop

\section{Embedding and first central limit result}
\label{sec:embedding}
In this section, a central limit result is established for the quantity
\begin{equation}
\label{eq:notyetteststatistic}
n^{5/6}\int_0^a\left\{\H_n(x)-\hat H_n(x)-\mu_n\right\}\,dF_0(x)=n^{5/6}\int_0^a\left\{\H_n(x-)-\hat H_n(x)-\mu_n\right\}\,dF_0(x),
\end{equation}
where $\mu_n$ denotes a centering sequence to be specified below. This result is the first step to be taken in order to obtain the limit result for $T_n$ defined in (\ref{eq:defTn}) (where the integral is taken with respect to $d\F_n$ rather than $dF_0$). In order to derive the limiting distribution of (\ref{eq:notyetteststatistic}),
we first replace the process $\H_n(x)-\hat H_n(x)$ by
$$
x\mapsto H_0(x)+\frac{E_n(x)}{\sqrt{n}\{1-F_0(x)}-\tilde H_n(x),\,x\in[0,a].
$$
where $E_n$ is the empirical process $\sqrt{n}\{\F_n-F_0\}$ and $\tilde H_n$ is the greatest convex minorant of the process
\begin{equation}
\label{eq:first_reduction}
x\mapsto H_0(x)+n^{-1/2}E_n(x)/\{1-F_0(x)\}\,,\,x\in[0,a].
\end{equation}
Next we use the strong approximation of the empirical process by a Brownian bridge $B_n$, yielding the approximation
$$
x\mapsto H_0(x)+\frac{n^{-1/2}B_n(F_0(x))}{1-F_0(x)},\,x\in[0,a],
$$
the process (\ref{eq:first_reduction}). This process is distributed as
\begin{equation}
\label{eq:first_BM}
x\mapsto H_0(x)+n^{-1/2}W\left(\frac{F_0(x)}{1-F_0(x)}\right),\,x\in[0,a],
\end{equation}
where $W$ is standard Brownian motion on $[0,\infty)$.
Next, the interval $[0,a]$ is split up in so-called big blocks separated by small blocks. The local contributions to the integral over the big blocks can be treated using the results of section \ref{sec:asproc}.

The first lemma to be proved states a contraction property for convex minorants that will be used repeatedly in the sequel. It is related to Marshall's Lemma in the theory of isotonic regression.
\begin{lemma}
\label{lemma:marshall2}
Let $f$ and $g$ be bounded functions on an interval $I\subset\R$ and let $C_f$ and $C_g$ be their greatest convex minorants, respectively. Then:
$$
\sup_{x\in I}\left|C_f(x)-C_g(x)\right|\le \sup_{x\in I}|f(x)-g(x)|.
$$
\end{lemma}

\noindent
{\bf Proof.} Using that $f\ge g-\sup_{u\in I}|f(u)-g(u)|$ and that $g\ge C_g$ by definition, it follows that $f\ge C_g-\sup_{u\in I}|f(u)-g(u)|$. Since the right hand side is convex, this means that it is {\it a} convex minorant of $f$ on $I$. Hence, it lies below the greatest convex minorant $C_f$ of $f$ on $I$:
$$
C_f(x)\ge C_g(x)-\sup_{u\in I}|f(u)-g(u)|,\,x\in I.
$$
Since this inequality also holds with $f$ and $g$ interchanged,
the result follows.\eop

We now consider the functional
\begin{align}
\label{functional_emp1}
&\int_{[0,a]}\left\{\H_n(x)-\hat H_n(x)\right\}\,dF_0(x)\nonumber\\
&=\int_{[0,a]}\left\{H_0(x)-\log\left(1-\frac{E_n(x)}{\sqrt{n}\bigl\{1-F_0(x)\bigr\}}\right)-\hat H_n(x)\right\}\,dF_0(x),
\end{align}
where $E_n=\sqrt{n}\{\F_n-F_0\}$ is the empirical process. The following lemma enables us to dispense with the logarithms.

\begin{lemma}
\label{lemma:logs-exit}
Let $\tilde H_n$ be the greatest convex minorant of the process
$$
x\mapsto H_0(x)+\frac{E_n(x)}{\sqrt{n}\bigl\{1-F_0(x)\bigr\}},\,x\in[0,a],
$$
where $F_0(a)<1$. Then:
\begin{enumerate}
\item[(i)]
$$
\int_{[0,a]}\left|\H_n(x)-H_0(x)-\frac{E_n(x)}{\sqrt{n}\bigl\{1-F_0(x)\bigr\}}\right|\,dF_0(x)=O_p\left(n^{-1}\right).
$$
\item[(ii)]
$$
\int_{[0,a]}\left|\hat H_n(x)-\tilde H_n(x)\right|\,dF_0(x)=O_p\left(n^{-1}\right).
$$
\end{enumerate}
\end{lemma}

\noindent
{\bf Proof.} (i). Let $A_n$ denote the event
$$
\left|\sup_{x\in[0,a]}\frac{E_n(x)}{\sqrt{n}\bigl\{1-F_0(x)\bigr\}}\right|\le\frac12.
$$
Then, by a well-known result in large deviation theory (``Chernoff's theorem"), we have
$$
\P\left(A_n^c\right)=O\left(e^{-nc}\right),
$$
for a constant $c>0$. If $A_n$ occurs, we can expand the logarithm, which yields:
$$
-\log\left\{1-\frac{E_n(x)}{\sqrt{n}\bigl\{1-F_0(x)\bigr\}}\right\}
=\frac{E_n(x)}{\sqrt{n}\bigl\{1-F_0(x)\bigr\}}+n^{-1}O\left(\sup_{x\in[0,a]}|E_n(x)|\right),
$$
and (i) now follows.\\
(ii). This follows from Lemma \ref{lemma:marshall2} and the argument of the proof of (i).\eop

\vspace{0.3cm}
We shall prove below that
$$
n^{5/6}\int_{[0,a]}\left\{H_0(x)+\frac{E_n(x)}{\sqrt{n}\bigl\{1-F_0(x)\bigr\}}-\tilde H_n(x)
-E\biggl\{H_0(x)+\frac{E_n(x)}{\sqrt{n}\bigl\{1-F_0(x)\bigr\}}-\tilde H_n(x)\biggr\}\right\}\,dF_0(x)
$$
converges in distribution to a normal distribution, which, together with Lemma \ref{lemma:logs-exit}, implies that
$$
n^{5/6}\int_{[0,a]}\left\{\H_n(x)-\hat H_n(x)
-E\biggl\{H_0(x)+\frac{E_n(x)}{\sqrt{n}\bigl\{1-F_0(x)\bigr\}}-\tilde H_n(x)\biggr\}\right\}\,dF_0(x)
$$
converges to the same normal distribution.

\begin{remark}
{\rm We avoid taking the expectation of
$$
\H_n(x)-\hat H_n(x),
$$
since $\H_n$ is infinite with a positive (but vanishing) probability on $[0,a]$, as is $\hat{H}_n$. This happens when the empirical distribution function $\F_n$ reaches the value $1$ on $[0,a]$.
}
\end{remark}

By Theorem 3 of \mycite{komlos:75} we can construct Brownian bridges $B_n$ on the same sample space as $\F_n$ such that
$$
Y_n=\sup_{x\in[0,a]}\frac{n^{1/2}\left|E_n(x)-B_n(F_0(x))\right|}{2\vee\log n}
$$
is a  random variable with with $EY_n\le C<\infty$ for all $n$.
Hence, for $n\ge 2$:
\begin{align}
\label{distance_emp_BB}
0\le E\sup_{x\in[0,a]}n^{-1/2}\left|\frac{E_n(x)}{1-F_0(x)}-\frac{B_n(F_0(x))}{1-F_0(x)}\right|\le \frac{EY_n\log n}{n(1-F_0(a))}
=O\left(\frac{\log n}{n}\right).
\end{align}
We now have the following result.

\begin{lemma}
\label{embedding_lemma}
Let $\tilde E_n$ be defined by
\begin{equation}
\label{tilde_E_n}
\tilde E_n(x)=\frac{B_n(F_0(x))}{1-F_0(x)}\,,\,x\in[0,a],
\end{equation}
and let $C_n^B$ be the greatest convex minorant of
$$
H_0(x)+n^{-1/2}\tilde E_n(x),\,x\in[0,a].
$$
Then
\begin{align}
\label{switch_to_BB}
&\int_{[0,a]}\left\{\H_n(x)-\hat H_n(x)\right\}\,dF_0(x)\nonumber\\
&=\int_{[0,a]}\left\{H_0(x)+n^{-1/2}\tilde E_n(x)-C_n^B(x)\right\}\,dF_0(x)+O_p\left(\frac{\log n}{n}\right).
\end{align}
\end{lemma}

\noindent
{\bf Proof.}
The result immediately follows from (\ref{distance_emp_BB}) and Lemmas \ref{lemma:marshall2} and \ref{lemma:logs-exit}(i).\eop

\vspace{0.3cm}
We now note that the process
$$
x\mapsto \frac{B(F_0(x))}{1-F_0(x)}\,,\,x\in[0,a],
$$
has the same distribution as the process
\begin{equation}
\label{V_n}
x\mapsto V_n(x)\stackrel{\mbox{\rm\small def}}=H_0(x)+n^{-1/2}W\left(\frac{F_0(x)}{1-F_0(x)}\right),\,x\in[0,a],
\end{equation}
where $W$ is standard Brownian motion on $\R_+$. So, if $C_n$ is the greatest convex minorant of the process
$$
x\mapsto H_0(x)+n^{-1/2}W\left(\frac{F_0(x)}{1-F_0(x)}\right)\,,\,x\in[0,a],
$$
we have:
\begin{align}
\label{switch_to_W}
&\int_{[0,a]}\left\{H_0(x)+n^{-1/2}\tilde E_n(x)-C_n^B(x)\right\}\,dF_0(x)\nonumber\\
&\stackrel{{\cal D}}=\int_{[0,a]}\left\{H_0(x)+n^{-1/2}W\left(\frac{F_0(x)}{1-F_0(x)}\right)-C_n(x)\right\}\,dF_0(x).
\end{align}

\begin{theorem}
\label{th:BB_limit}
Let $h_0$ be strictly positive on $[0,a]$, with a strictly positive continuous derivative $h_0'$ on $(0,a)$,
 which also has a strictly positive right limit at $0$ and a strictly positive left limit at $a$.
Moreover, let $S_n$ be defined by
\begin{equation}
\label{def_S_n}
S_n(x)=H_0(x)+n^{-1/2}W\left(\frac{F_0(x)}{1-F_0(x)}\right)-C_n(x),\,x\in[0,a],
\end{equation}
where $C_n$ is the greatest convex minorant of $V_n$ defined in (\ref{V_n})
and let $D_n$ be defined by
\begin{equation}
\label{def_D_n}
D_n=\int_0^a S_n(x)\,dF_0(x).
\end{equation}
Finally, let $C(0)$ and $\s^2$ be defined as in Theorem \ref{th:BM-limit}. Then:
\begin{align*}
n^{5/6}\left\{D_n-ED_n\right\}\stackrel{{\cal D}}\longrightarrow N(0,\s_{H_0}^2),\,n\to\infty,
\end{align*}
where
\begin{equation}
\label{asymp_expectation1}
n^{2/3}ED_n\to E|C(0)|\int_0^a\left(\frac{2h_0(t)f_0(t)}{h_0'(t)}\right)^{1/3}\,dH_0(t),
\end{equation}
and
\begin{equation}
\label{asymp_var1}
\s_{H_0}^2=2^{4/3}\s^2\int_0^{a}\frac{h_0(t)^2\{h_0(t)f_0(t)\}^{1/3}}{h_0'(t)^{4/3}}\,dH_0(t)
\end{equation}
\end{theorem}

The following corollary is immediate from the preceding.

\begin{corollary}
Let $h_0$ be strictly positive on $[0,a]$, with a strictly positive continuous derivative $h_0'$ on $(0,a)$,
 which also has a strictly positive right limit at $0$ and a strictly positive left limit at $a$.
Then:
\begin{align*}
n^{5/6}\left\{\int_0^a\left\{\H_n(x)-\hat H_n(x)\right\}\,dF_0(x)-ED_n\right\}\stackrel{{\cal D}}\longrightarrow N(0,\s_{H_0}^2),\,n\to\infty,
\end{align*}
where $ED_n$ and $\s_{H_0}^2$ are defined as in Theorem \ref{th:BB_limit}.
\end{corollary}

For the proof of Theorem \ref{th:BB_limit} we divide the interval $[0,a]$ into $m_n$  intervals $I_{n,k}$ with (equal) length of order $n^{-1/3}\log n$ (big blocks), separated by intervals $J_{n,k}$ ($k=2,3,\ldots,m_n$) with length of order $2n^{-1/3}\sqrt{\log n}$ (small blocks). The small interval $J_{n,1}$ to the left of $I_{n,1}$  has half the length of the other separating blocks as has the small interval $J_{n,m_n+1}$ to the right of $I_{n,m_n}$. Hence,
 $$
[0,a]=J_{n,1}\cup I_{n,1}\cup J_{n,2}\cup I_{n,2}\cdots\cup J_{n,m_n}\cup I_{n,m_n}\cup J_{n,m_n+1}.
$$
For $k=2,3,\ldots,m_n$, let $\tilde J_{n,k}$ be the interval with the same right endpoint as $J_{n,k}$ with half the length of $J_{n,k}$ and take $\tilde{J}_{n,1}=J_{n,1}$.  For $k=1,2,\ldots,m_n-1$ let  $\bar J_{n,k+1}$ be the interval with the same left endpoint as $J_{n,k+1}$ with half the length of $J_{n,k+1}$ and $\tilde{J}_{n,m_n+1}=J_{n,m_n+1}$. Then
$$
[0,a]= \tilde{J}_{n,1}\cup  I_{n,1}\cup \bar{J}_{n,2}\cup \tilde{J}_{n,2}\cup I_{n,2}\cdots\cup \tilde{J}_{n,m_n}\cup I_{n,m_n}\cup \bar{J}_{n,m_n+1}
$$
where all $I$-intervals have the same length, of order $n^{-1/3}\log n$ and the $J$-intervals have the same length of (smaller) order $n^{-1/3}\sqrt{\log n}$.
Finally, let the interval $L_{n,k}$ be defined by
\begin{equation}
\label{eq:intervals}
L_{n,k}=\tilde J_{n,k}\cup I_{n,k}\cup \bar J_{n,k+1}=[a_{nk},a_{n,k+1}),\,\,\, k=1,2,\ldots,m_n, \mbox{ yielding } [0,a)=\cup_{k=1}^{m_n}L_{n,k}.
\end{equation}
Note that $m_n\sim a n^{1/3}/\log{n}$ and see the figure below for the structure of the partition.

\vspace{0.5cm}
\begin{picture}(100,20)
\put(0,-2){\line(0,1){4}}
\put(0,0){\line(1,0){200}}
\multiput(200,0)(4,0){13}{\line(1,0){2}}
\put(240,0){\line(1,0){180}}
\put(0,-15){$0$}
\put(1,10){$\tilde{J}_{n,1}$}
\put(20,-2){\line(0,1){4}}
\put(36,10){$I_{n,1}$}
\put(70,-2){\line(0,1){4}}
\put(71,10){$\bar{J}_{n,2}$}
\put(90,-10){\line(0,1){12}}
\put(91,10){$\tilde{J}_{n,2}$}
\put(110,-2){\line(0,1){4}}
\put(126,10){$I_{n,2}$}
\put(140,-10){\vector(1,0){40}}
\put(140,-10){\vector(-1,0){50}}
\put(126,-20){$L_{n,2}$}
\put(160,-2){\line(0,1){4}}
\put(161,10){$\bar{J}_{n,3}$}
\put(180,-10){\line(0,1){12}}
\put(420,-15){$a$}
\put(420,-2){\line(0,1){4}}
\put(395,10){$\bar{J}_{n,m_n+1}$}
\put(400,-2){\line(0,1){4}}
\put(360,10){$I_{n,m_n}$}
\put(350,-2){\line(0,1){4}}
\put(330,-2){\line(0,1){4}}
\put(310,-2){\line(0,1){4}}
Ê\end{picture}

\vspace{1cm}

The (key) localization lemma below  and proved in the appendix, shows that on intervals $I_{n,k}$ the global convex minorant of $V_n$ (defined in (\ref{V_n})) over $[0,a]$ coincides with high probability with the restriction to $I_{nk}$ of the local convex minorant of the process $V_n$ on the interval $L_{n,k}$.

\begin{lemma}
\label{localization_lemma}
Let $h_0$ be strictly positive on $[0,a]$, with a strictly positive continuous derivative $h_0'$ on $(0,a)$,
 which also has a strictly positive right limit at $0$ and a strictly positive left limit at $a$.
Then:
\begin{enumerate}
\item[(i)]
The probability that there exists a $k$, $1\le k\le m_n$, such that  the greatest convex minorant $C_n$ of $V_n$ is different on the interval $I_{nk}$ from the restriction to $I_{nk}$ of the (local) greatest convex minorant of $V_n$ on $L_{nk}$, is bounded above by
$$
c_1\exp\left\{-c_2(\log n)^{3/2}\right\},
$$
for constants $c_1,c_2>0$, uniformly in $n$.
\item[(ii)] The probability that there exists a $k$, $1\le k\le m_n$, such that   $C_n$ has no change of slope in an interval $\bar J_{nk}$ or $\tilde J_{nk}$ is also bounded by
$$
c_1\exp\left\{-c_2(\log n)^{3/2}\right\},
$$
for constants $c_1,c_2$, uniformly in $n$.
\end{enumerate}
\end{lemma}

For each $n\ge 1$ and $1\le k\le m_n$ define independent standard Brownian motions $W_{n1},\dots,W_{n,m_n}$ and consider the processes
$$
x\mapsto H_0(x)-H_0(a_{nk})+n^{-1/2}W_{nk}\left(\frac{F_0(x)}{1-F_0(x)}-\frac{F_0(a_{nk})}{1-F_0(a_{nk})}\right),\,x\in L_{nk}.
$$
Denote the greatest convex minorants of these processes (on $L_{nk}$)  by $C_{nk}$.
Furthermore,  define the processes $S_{nk}$ by
\begin{equation}
\label{S_{nk}}
S_{nk}(x)=H_0(x)-H_0(a_{nk})+n^{-1/2}W_{nk}\left(\frac{F_0(x)}{1-F_0(x)}-\frac{F_0(a_{nk})}{1-F_0(a_{nk})}\right)-C_{nk}(x),\,x\in L_{nk}.
\end{equation}

\begin{lemma}
\label{lemma:central_limit1}
Assume that the conditions of Theorem \ref{th:BB_limit} are satisfied. Moreover, let $C(0)$  be defined as in Theorem \ref{th:BM-limit} and  $\s_{H_0}^2$  as in Theorem \ref{th:BB_limit}. Then:
\begin{align*}
&n^{5/6}\sum_{k=1}^{m_n}\int_{I_{nk}}\left\{S_{nk}(x)-E S_{nk}(x)\right\}\,dF_0(x)
\stackrel{{\cal D}}\longrightarrow N(0,\s_{H_0}^2),\,n\to\infty,
\end{align*}
where (see (\ref{asymp_expectation1}))
\begin{equation}
\label{sum_local_expectations}
n^{2/3}\sum_{k=1}^{m_n}\int_{I_{nk}}E S_{nk}(x)\,dF_0(x)
\to E|C(0)|\int_0^a\left(\frac{2 h_0(t)f_0(t)}{h_0'(t)}\right)^{1/3}\,dH_0(t)\,,\,n\to\infty,
\end{equation}
\end{lemma}

\noindent
{\bf Proof.}
Let $c_n=n^{1/3}|L_{nk}|\sim\log n$ and $I_{nk}=[a_{nk}+n^{-1/3}\sqrt{c_n},a_{nk}+n^{-1/3}(c_n-\sqrt{c_n})]$. We then have:
\begin{align*}
&n\int_{I_{nk}}\left\{S_{nk}(x)-ES_{nk}(x)\right\}\,dF_0(x)\\
&=\int_{\sqrt{c_n}}^{c_n-\sqrt{c_n}}\left\{n^{1/6}W_{nk}\left(\frac{F_0(a_{nk}+n^{-1/3}x)}{1-F_0(a_{nk}+n^{-1/3}x)}-\frac{F_0(a_{nk})}{1-F_0(a_{nk})}\right)-n^{2/3}C_{nk}(a_{nk}+n^{-1/3}x)\right.\\
&\qquad\qquad\left.-E\left\{n^{1/6}W_{nk}\left(\frac{F_0(a_{nk}+n^{-1/3}x)}{1-F_0(a_{nk}+n^{-1/3}x)}-\frac{F_0(a_{nk})}{1-F_0(a_{nk})}\right)-n^{2/3}C_{nk}(a_{nk}+n^{-1/3}x)\right\}\right\}\\
&\qquad\qquad\qquad\qquad\qquad\qquad\qquad\qquad\qquad\qquad\qquad\qquad\qquad\cdot\,f_0\left(a_{nk}+n^{-1/3}x\right)\,dx,
\end{align*}
Here we use that the (first two) deterministic terms in (\ref{S_{nk}}) drop out because of subtraction of the expectation. This implies that
\begin{align*}
&n\int_{I_{nk}}\left\{S_{nk}(x)-ES_{nk}(x)\right\}\,dF_0(x)\\
&\stackrel{{\cal D}}=\int_{\sqrt{c_n}}^{c_n-\sqrt{c_n}}\left\{n^{1/6}W\left(\frac{F_0(a_{nk}+n^{-1/3}x)}{1-F_0(a_{nk}+n^{-1/3}x)}-\frac{F_0(a_{nk})}{1-F_0(a_{nk})}\right)-C_{nk}'(x)\right.\\
&\qquad\left.-E\left\{n^{1/6}W\left(\frac{F_0(a_{nk}+n^{-1/3}x)}{1-F_0(a_{nk}+n^{-1/3}x)}-\frac{F_0(a_{nk})}{1-F_0(a_{nk})}\right)-C_{nk}'(x)\right\}\right\}
\cdot\,f_0\left(a_{nk}+n^{-1/3}x\right)\,dx,
\end{align*}
where $C_{nk}'$ is the greatest convex minorant of the process
\begin{align*}
x&\mapsto n^{2/3}\left\{H_0(a_{nk}+n^{-1/3}x)-H_0(a_{nk})-n^{-1/3}xh_0(a_{nk})\right\}\\
&\qquad\qquad\qquad+n^{1/6}W\left(\frac{F_0(a_{nk}+n^{-1/3}x)}{1-F_0(a_{nk}+n^{-1/3}x)}-\frac{F_0(a_{nk})}{1-F_0(a_{nk})}\right),\,\, x\in[0,c_n].
\end{align*}
Here we use that adding a linear function to a function does not change the difference between this function and its greatest convex minorant.
Note that the integrals on $I_{nk}$ only depend on the increments of the Brownian motion process on the corresponding disjoint intervals $L_{nk}$ and therefore are independent.
For the individual integrals we are close to the situation of Theorem \ref{th:BM-limit3}, with, for $c_n\rightarrow\infty$, on $[0,c_n]$ (note that $n$ is determined by $c_n$,  $n=e^{c_n}$)
\begin{align*}
&F_{c_n}(x)=n^{1/3}\left\{F_0(a_{nk}+n^{-1/3}x)-F_0(a_{nk})\right\}=f_0(a_{nk})x(1+o(1)),\\ &H_{c_n}(x)=n^{2/3}\left\{H_0(a_{nk}+n^{-1/3}x)-H_0(a_{nk})-n^{-1/3}xh_0(a_{nk})\right\}=\tfrac12h_0^\prime(a_{nk})x^2(1+o(1))\\
&\mbox{ and } G_{c_n}(x)=n^{1/3}\left\{\frac{F_0(a_{nk}+n^{-1/3}x)}{1-F_0(a_{nk}+n^{-1/3}x)}-\frac{F_0(a_{nk})}{1-F_0(a_{nk})}\right\}=\frac{f_0(a_{nk})x}{(1-F_0(a_{nk}))^2}(1+o(1))
\end{align*}
This yields:
$$
\mbox{var}\left(\frac{n}{\sqrt{c_n}}\int_{I_{nk}}S_{nk}(x)\,dF_0(x)\right)
\sim\s_{nk}^2,\,n\to\infty,
$$
uniformly in $k=1,\dots,m_n$, where
\begin{align*}
\s_{nk}^2&=\frac{\left(f_0(a_{nk})/\{1-F_0(a_{nk})\}^2\right)^{5/3}f_0(a_{nk})^2}{\left(\tfrac12h_0'(a_{nk})\right)^{4/3}}\sigma^2\\
&=\frac{\left(h_0(a_{nk})\right)^{10/3}f_0(a_{nk})^{1/3}}{\left(\tfrac12h_0'(a_{nk})\right)^{4/3}}
\sigma^2=\frac{2^{4/3}h_0(a_{nk})^3\left\{h_0(a_{nk})f_0(a_{nk})\right\}^{1/3}}{h_0'(a_{nk})^{4/3}}
\sigma^2,
\end{align*}
and  $\sigma^2$ is defined as in Theorem \ref{th:BM-limit}. Likewise, also with $C(0)$ as defined in Theorem \ref{th:BM-limit},
\begin{align*}
\frac{n}{c_n}\int_{I_{nk}} E S_{nk}(x)\,dF_0(x)\sim
\frac{2^{1/3}h_0(a_{nk})\left\{h_0(a_{nk})f_0(a_{nk})\right\}^{1/3}E|C(0)|}{h_0'(a_{nk})^{1/3}}
\end{align*}

Since the fourth moments of
$$
\frac{n}{\sqrt{c_n}}\int_{I_{nk}}\left\{S_{nk}(x)-ES_{nk}(x)\right\}\,dF_0(x)
$$
are uniformly bounded by Theorem \ref{th:BM-limit3}, we get for each $\e>0$, using Chebyshev's inequality,
\begin{align*}
&\sum_{k=1}^{m_n}\P\left\{m_n^{-1/2}\left|\frac{n}{\sqrt{c_n}}\int_{I_{nk}}\left\{S_{nk}(x)-ES_{nk}(x)\right\}\,dF_0(x)\right|\ge \e\right\}\to0,\,n\to\infty.
\end{align*}
Using that $m_n^{-1}\sim a^{-1} n^{-1/3}\log n$
and that the intervals $I_{nk}$ have lengths of order $n^{-1/3}\log n$, we get:
\begin{align*}
m_n^{-1}\sum_{k=1}^{m_n}\s_{nk}^2
&\sim m_n^{-1}\sum_{k=1}^{m_n}\frac{2^{4/3}h_0(a_{nk})^3\left\{h_0(a_{nk})f_0(a_{nk})\right\}^{1/3}}{h_0'(a_{nk})^{4/3}}\sigma^2\\
&\longrightarrow \frac{2^{4/3}\sigma^2}{a}\int_0^{a}\frac{h_0(t)^3\{h_0(t)f_0(t)\}^{1/3}}{h_0'(t)^{4/3}}\,dt.
\end{align*}
Since $m_n=an^{1/3}/c_n$, the normal convergence criterion on p.\ 316 of \mycite{loeve:63} now gives:
\begin{align*}
&n^{5/6}\sum_{k=1}^{m_n}\int_{I_{nk}}\left\{S_{nk}(x)-ES_{nk}(x)\right\}\,dF_0(x)\\
&=m_n^{-1/2}\sum_{k=1}^{m_n}\frac{n\sqrt{a}}{\sqrt{c_n}}\int_{I_{nk}}\left\{S_{nk}(x)-ES_{nk}(x)\right\}\,dF_0(x)\stackrel{{\cal D}}\longrightarrow N\left(0,\s^2_{H_0}\right).
\end{align*}
Also note that:
\begin{align*}
&m_n^{-1/2}\sum_{k=1}^{m_n}\frac{n}{c_n^{1/2}}\int_{I_{nk}} E S_{nk}(x)\,dF_0(x)
\sim m_n^{-1/2}c_n^{1/2}\sum_{k=1}^{m_n}\frac{2^{1/3}h_0(a_{nk})\left\{h_0(a_{nk})f_0(a_{nk})\right\}^{1/3}}{h_0'(a_{nk})^{1/3}}\\
&\sim \sqrt{m_nc_n}E|C(0)|\int_0^a\frac{2^{1/3}h_0(t)\left\{h_0(t)f_0(t)\right\}^{1/3}}{h_0'(t)^{1/3}}\,dt  =n^{1/6}\int_0^a\left(\frac{2h_0(t)f_0(t)}{h_0'(t)}\right)^{1/3}\,dH_0(t).
\end{align*}
\eop

In applications of this result, used in a bootstrap approach to the computation of critical values, we need the following lemma, which gives a more precise expansion of the asymptotic representation of the expectation, given in (\ref{sum_local_expectations}).

\begin{lemma}
\label{lemma:remainder term}
Assume that the conditions of Theorem \ref{th:BB_limit} are satisfied and assume in addition that $h_0$ has a bounded second derivative on $[0,a]$. Then
\begin{equation}
\label{sum_local_expectations_refined}
n^{2/3}ED_n
=E|C(0)|\int_0^a\left(\frac{2 h_0(t)f_0(t)}{h_0'(t)}\right)^{1/3}\,dH_0(t)+o\left(n^{-1/6}\right).
\end{equation}
\end{lemma}

\noindent
{\bf Proof.}  We have:
\begin{align}
\label{rem_terms}
\frac{n}{c_n}\int_{I_{nk}} ES_{nk}(x)\,dF_0(x)
=\frac{h_0(a_{nk})\bigl\{2h_0(a_{nk})f_0(a_{nk})\bigr\}^{1/3}E|\zeta(0)|}{h_0'(a_{nk})^{1/3}}+o\left(n^{-1/6}\right),
\end{align}
uniformly in $k=1,\dots,m_n$. This is seen in the following way.

On the intervals $I_{nk}$ we get:
\begin{align*}
F_{c_n}(x)&\stackrel{\mbox{\small def}}=n^{1/3}\left\{F_0(a_{nk}+n^{-1/3}x)-F_0(a_{nk})\right\}=f_0(a_{nk})x+O\left(n^{-1/3}(\log n)^2\right),\\
H_{c_n}(x)&\stackrel{\mbox{\small def}}=n^{2/3}\left\{H_0(a_{nk}+n^{-1/3}x)-H_0(a_{nk})-n^{-1/3}xh_0(a_{nk})\right\}\\
&=\tfrac12h_n^\prime(a_{nk})x^2+O\left(n^{-1/3}(\log n)^3\right),\\
G_{c_n}(x)&\stackrel{\mbox{\small def}}=n^{1/3}\left\{\frac{F_0(a_{nk}+n^{-1/3}x)}{1-F_0(a_{nk}+n^{-1/3}x)}-\frac{F_0(a_{nk})}{1-F_0(a_{nk})}\right\}
=\frac{f_0(a_{nk})x}{(1-F_0(a_{nk}))^2}+O\left(n^{-1/3}(\log n)^2\right).
\end{align*}
uniformly in $k=1,\dots,m_n$.
The relation for $F_{c_n}$ and $G_{c_n}$ immediately follow from the mean value theorem, applied on the remainder term, together with the conditions of Theorem \ref{th:BB_limit}, which yield that $h_0'$ and $f_0'$ are uniformly bounded. In the expansion of $H_{c_n}$ we use the boundedness of the second derivative $h_0''$. 

Combining these relations gives (\ref{rem_terms}), and hence:
\begin{align*}
&n^{2/3}\sum_{k=1}^{m_n}\int_{I_{nk}} ES_{nk}(x)\,dF_0(x)
=\frac{n}{c_n}\sum_{k=1}^{m_n}\int_{I_{nk}} ES_{nk}(x)\,dF_0(x)\left|I_{nk}\right|\\
&=\frac{n}{c_n}\sum_{k=1}^{m_n}\frac{h_0(a_{nk})\bigl\{2h_0(a_{nk})f_0(a_{nk})\bigr\}^{1/3}E|\zeta(0)|}{h_0'(a_{nk})^{1/3}}
\left|I_{nk}\right|+O\left(n^{-1/3}(\log n)^3\sum_{k=1}^{m_n}\left|I_{nk}\right|\right)\\
&=E|C(0)|\sum_{k=1}^{m_n}\int_{I_{nk}}\left(\frac{2 h_0(t)f_0(t)}{h_0'(t)}\right)^{1/3}\,dH_0(t)+O\left(n^{-1/3}(\log n)^3\right)
\end{align*}
where, in the last line, we use again the boundedness of $f_0'$, $h_0'$ and $h_0''$, combined with the mean value theorem on the intervals $I_{nk}$. Note that part (i) of Lemma \ref{localization_lemma} tells us that the probability that $S_n$ is different from $S_{nk}$ on $I_{nk}$ is bounded above by
$$
c_1\exp\left\{-c_2(\log n)^{3/2}\right\},
$$
so we also have:
\begin{align*}
&n^{2/3}\sum_{k=1}^{m_n}\int_{I_{nk}} ES_n(x)\,dF_0(x)
&=E|C(0)|\sum_{k=1}^{m_n}\int_{I_{nk}}\left(\frac{2 h_0(t)f_0(t)}{h_0'(t)}\right)^{1/3}\,dH_0(t)+O\left(n^{-1/3}(\log n)^3\right).
\end{align*}
It is clear that we get in a similar way:
\begin{align*}
&n^{2/3}\sum_{k=1}^{m_n}\int_{L_{nk}\setminus I_{nk}} ES_n(x)\,dF_0(x)\\
&=E|C(0)|\sum_{k=1}^{m_n}\int_{L_{nk}\setminus I_{nk}}\left(\frac{2 h_0(t)f_0(t)}{h_0'(t)}\right)^{1/3}\,dH_0(t)+O\left(n^{-1/3}(\log n)^{3/2}\right).
\end{align*}
Hence
\begin{align*}
n^{2/3}ED_n=n^{2/3}\sum_{k=1}^{m_n}\int_{L_{nk}} ES_n(x)\,dF_0(x)=
\int_0^a\left(\frac{2 h_0(t)f_0(t)}{h_0'(t)}\right)^{1/3}\,dH_0(t)+o\left(n^{-1/6}\right).
\end{align*}
\eop

\vspace{0.3cm}
We can now prove Theorem \ref{th:BB_limit}.

\vspace{0.3cm}
\noindent
{\bf Proof of Theorem \ref{th:BB_limit}.} By Lemmas \ref{localization_lemma} and \ref{lemma:central_limit1} we have:
$$
n^{5/6}\sum_{k=1}^{m_n}\int_{I_{nk}}\left\{S_{nk}(x)-E S_{nk}(x)\right\}\,dF_0(x)
\stackrel{{\cal D}}\longrightarrow N(0,\s_{H_0}^2),\,n\to\infty.
$$
By part (i) of Lemma \ref{localization_lemma}, the probability that $S_n$ is different from $S_{nk}$ on $I_{nk}$ is bounded above by
$$
c_1\exp\left\{-c_2(\log n)^{3/2}\right\},
$$
implying that also:
$$
n^{5/6}\sum_{k=1}^{m_n}\int_{I_{nk}}\left\{S_n(x)-E S_n(x)\right\}\,dF_0(x)
\stackrel{{\cal D}}\longrightarrow N(0,\s_{H_0}^2),\,n\to\infty.
$$
For similar reasons we have:
$$
n^{5/6}\sum_{k=1}^{m_n}\int_{L_{nk}\setminus I_{nk}}\left\{S_n(x)-E S_n(x)\right\}\,dF_0(x)
\stackrel{p}\longrightarrow 0,\,n\to\infty,
$$
where we use Theorem \ref{th:BM-limit3} (this is the essence of the ``big blocks, small blocks" method).
The result now follows, since
$$
D_n=n^{5/6}\sum_{k=1}^{m_n}\int_{L_{nk}}\left\{S_n(x)-E S_n(x)\right\}\,dF_0(x).
$$
\eop

\section{Further central limit results}
\label{sec:further}
In order to derive the asymptotic distribution of the statistic $U_n$ defined in (\ref{eq:defUn}) and used in the simulations in \mycite{GrJo10}, we first consider the statistic
$$
\int_{[0,a]}\bigl\{\F_n(x)-\hat F_n(x)\bigr\}\,dF_0(x),
$$
which is analogous to the statistic discussed in the preceding section, but
has $\F_n(x)-\hat F_n(x)$ as integrand instead of $\H_n(x)-\hat H_n(x)$.
We have, if $E_n$ again denotes the empirical process, and arguing as in the proof of Lemma \ref{lemma:logs-exit} (i),
\begin{align*}
&\F_n(x)=1-\exp\left\{-\H_n(x)\right\}=1-\exp\left\{-H_0(x)+\log\left\{1-\frac{n^{-1/2}E_n(x)}{1-F_0(x)}\right\}\right\}\\
&=1-\exp\left\{-H_0(x)-\frac{n^{-1/2}E_n(x)}{1-F_0(x)}\right\}+O_p\left(n^{-1}\right),
\end{align*}
uniformly for $x\in[0,a]$. Hence, defining, as in Lemma \ref{lemma:logs-exit}, $H_n$ as the greatest convex minorant of the process
$$
x\mapsto H_0(x)+\frac{n^{-1/2}E_n(x)}{1-F_0(x)}\,,
$$
we get, by Lemma \ref{lemma:logs-exit},
\begin{align*}
&\F_n(x)-\hat F_n(x)=\exp\left\{-H_n(x)\right\}
-\exp\left\{-H_0(x)-\frac{n^{-1/2}E_n(x)}{1-F_0(x)}\right\}+O_p\left(n^{-1}\right)\\
&=\exp\left\{-H_n(x)\right\}\left\{1-\exp\left\{-H_0(x)-\frac{n^{-1/2}E_n(x)}{1-F_0(x)}+H_n(x)\right\}\right\}+O_p\left(n^{-1}\right).
\end{align*}
Next, replacing $E_n(x)$ by $B_n(F_0(x))$, as in Lemma \ref{embedding_lemma}, where $(B_n)$ are the approximating Brownian bridges, we get:
\begin{align*}
&\F_n(x)-\hat F_n(x)\nonumber\\
&=\exp\left\{-C_n^B(x)\right\}\left\{1-\exp\left\{-H_0(x)-\frac{n^{-1/2}B_n(F_0(x))}{1-F_0(x)}+C_n^B(x)\right\}\right\}+O_p\left(\frac{\log n}{n}\right),
\end{align*}
where $C_n^B$ is the greatest convex minorant of the process
$$
x\mapsto H_0(x)+\frac{n^{-1/2}B_n(F_0(x))}{1-F_0(x)}\,,\,x\in[0,a].
$$
Again using the results of the preceding section, it is seen that this implies that
\begin{align}
\label{BM_replacement}
&\F_n(x)-\hat F_n(x)\nonumber\\
&\stackrel{{\cal D}}
=\exp\left\{-C_n(x)\right\}\left\{1-\exp\left\{-H_0(x)-n^{-1/2}W\left(\frac{F_0(x)}{1-F_0(x)}\right)+C_n(x)\right\}\right\}+O_p\left(\frac{\log n}{n}\right),
\end{align}
where $W$ is standard Brownian motion on $[0,\infty)$, and $C_n$ is the greatest convex minorant of the process
$$
x\mapsto H_0(x)+n^{-1/2}W\left(\frac{F_0(x)}{1-F_0(x)}\right)\,,\,x\in[0,a].
$$
This representation suggests to consider
\begin{align*}
&\exp\left\{-H_0(x)\right\}\left\{H_0(x)+n^{-1/2}W\left(\frac{F_0(x)}{1-F_0(x)}\right)-C_n(x)\right\}\\
&=\left\{1-F_0(x)\right\}\left\{H_0(x)+n^{-1/2}W\left(\frac{F_0(x)}{1-F_0(x)}\right)-C_n(x)\right\},\,x\in[0,a].
\end{align*}
We have the following result.

\begin{lemma}
\label{lemma:df_limit}
Let $h_0$ be strictly positive on $[0,a]$, with a strictly positive continuous derivative $h_0'$ on $(0,a)$,
which also has a strictly positive right limit at $0$ and a strictly positive left limit at $a$.
Moreover, let $C_n, S_n$ and $V_n$ be defined as in Theorem \ref{th:BB_limit}
and let $D_n^{F_0}$ be defined by
\begin{equation}
\label{def_D_n_F_0}
D_n^{F_0}=\int_0^a S_n(x)\bigl\{1-F_0(x)\bigr\}\,dF_0(x),
\end{equation}
Then:
\begin{align*}
n^{5/6}\left\{D_n^{F_0}-ED_n^{F_0}\right\}\stackrel{{\cal D}}\longrightarrow N(0,\s_{F_0}^2),\,n\to\infty,
\end{align*}
where
$$
ED_N^{F_0}\sim n^{-2/3}E|C(0)|\int_0^a\left(\frac{2h_0(t)f_0(t)}{h_0'(t)}\right)^{1/3}dF_0(t),\,n\to\infty,
$$
$$
\s_{F_0}^2=\s^2\int_0^{a}\left(\frac{2h_0(t)f_0(t)}{h_0'(t)}\right)^{4/3}dF_0(t),
$$
and $\s^2$ is defined as in Theorem \ref{th:BM-limit}.
\end{lemma}

\noindent
{\bf Proof.} The only difference with Theorem \ref{th:BB_limit} is that $dF_0(t)$ is replaced by $\{1-F_0(t)\}\,dF_0(t)$ in the integral. This means that instead of $ED_n$ we get:
$$
ED_N^{F_0}\sim n^{-2/3}E|C(0)|\int_0^a\frac{2^{1/3}h_0(t)\left\{h_0(t)f_0(t)\right\}^{1/3}\{1-F_0(t)\}}{h_0'(t)^{1/3}}\,dt,\,n\to\infty,
$$
and instead of $\s^2_{H_0}$ we get:
$$
\s_{F_0}^2=2^{4/3}\s^2\int_0^{a}\frac{h_0(t)^3\{h_0(t)f_0(t)\}^{1/3}}{h_0'(t)^{4/3}}\bigl\{1-F_0(t)\}^2\,dt
=\s^2\int_0^{a}\left(\frac{2h_0(t)f_0(t)}{h_0'(t)}\right)^{4/3}\,dF_0(t).
$$
\eop

\vspace{0.3cm}
We have the following corollary.

\begin{corollary}
\label{cor:df_limit}
Let $h_0$ be strictly positive on $[0,a]$, with a strictly positive continuous derivative $h_0'$ on $(0,a)$,
which also has a strictly positive right limit at $0$ and a strictly positive left limit at $a$.
Moreover, let $S_n'$ be defined by:
$$
S_n'(x)=\F_n(x)-\hat F_n(x),\,x\in[0,a],
$$
where $\hat F_n$ is defined in (\ref{eq:defUn}) and let $D_n'$ be defined by
\begin{equation}
\label{def_D_n'}
D_n'=\int_0^a S_n'(x)\,dF_0(x),
\end{equation}
Then,
\begin{align*}
n^{5/6}\left\{D_n'-ED_n'\right\}\stackrel{{\cal D}}\longrightarrow N(0,\s_{F_0}^2),\,n\to\infty,
\end{align*}
where $\s_{F_0}^2$ is defined as in Lemma \ref{lemma:df_limit}.
\end{corollary}

\noindent
{\bf Proof.} This is (in a sense) an application of the delta method. By (\ref{BM_replacement}) we can replace $\F_n-\hat F_n$ by:
$$
\exp\left\{-C_n(x)\right\}\left\{1-\exp\left\{-H_0(x)-n^{-1/2}W\left(\frac{F_0(x))}{1-F_0(x)}\right)+C_n(x)\right\}\right\}.
$$
We also have, using notation of the same type as in the proof of Lemma \ref{lemma:central_limit1},
\begin{align*}
&\int_0^a E\left\{H_0(x)-C_n(x)\right\}^2\,dF_0(x)\\
&\sim \sum_{k=1}^{m_n}\int_0^{c_n}E\left\{H_0(a_{nk}+n^{-1/3}u)-H_0(a_{nk})-C_n(a_{nk}+n^{-1/3}u)
+C_n(a_{nk})\right\}^2\,f_0(a_{nk})\,du\\
&\sim n^{-5/3}\sum_{k=1}^{m_n}\int_0^{c_n}E\left\{\tfrac12h_0'(a_{nk})u^2-C_{nk}(u)\right\}^2\,f_0(a_{nk})\,du,
\end{align*}
where $C_{nk}$ is the greatest convex minorant of the process
$$
x\mapsto \tfrac12h_0'(a_{nk})u^2+W\left(\frac{h_0(a_{nk})u}{1-F_0(a_{nk})}\right),\,u\in[0,c_n].
$$
By Brownian scaling, we get:
\begin{align*}
&\int_0^{c_n}E\left\{\tfrac12h_0'(a_{nk})u^2-C_{nk}(u)\right\}^2\,f_0(a_{nk})\,du\\
&\sim c_n f_0(a_{nk})\left(\tfrac12h_0'(a_{nk})\right)^{2/3}\left(\frac{h_0(a_{nk})}{1-F_0(a_{nk})}\right)^{4/3} EC(0)^2,
\end{align*}
where $C$ is the greatest convex minorant of $x\mapsto W(x)+x^2,\,x\in\R$. So we find:
\begin{align}
\label{CS1}
&\int_0^a E\left\{H_0(x)-C_n(x)\right\}^2\,dF_0(x)
\sim n^{-4/3}EC(0)^2\int_0^a\left(\tfrac12h_0'(t)\right)^{2/3}\left(\frac{h_0(t)}{1-F_0(t)}\right)^{4/3}\,dF_0(t).
\end{align}
We also have:
\begin{align}
\label{CS2}
&\int_0^a E\left\{H_0(x)+n^{-1/2}W\left(\frac{F_0(x)}{1-F_0(x)}\right)-C_n(x)\right\}^2\,dF_0(x)\nonumber\\
&\sim n^{-4/3}EC(0)^2\int_0^a\left(\tfrac12h_0'(t)\right)^{2/3}\left(\frac{h_0(t)}{1-F_0(t)}\right)^{4/3}\,dF_0(t).
\end{align}
Hence, by (\ref{CS1}) and (\ref{CS2}),
\begin{align*}
&\int_0^a\exp\left\{-C_n(x)\right\}\left\{1-\exp\left\{-H_0(x)-n^{-1/2}W\left(\frac{F_0(x))}{1-F_0(x)}\right)+C_n(x)\right\}\right\}\,dF_0(t)\\
&=\int_0^a\left\{1-F_0(t)\right\}\left\{1-\exp\left\{-H_0(x)-n^{-1/2}W\left(\frac{F_0(x))}{1-F_0(x)}\right)+C_n(x)\right\}\right\}\,dF_0(t)+O_p\left(n^{-4/3}\right)\\
&=\int_0^a\left\{1-F_0(t)\right\}\left\{H_0(x)+n^{-1/2}W\left(\frac{F_0(x))}{1-F_0(x)}\right)-C_n(x)\right\}\,dF_0(t)+O_p\left(n^{-4/3}\right),
\end{align*}
where we also use the Cauchy-Schwarz inequality in the first equality.

For the expectation we get similarly
\begin{align*}
&\int_0^a E\left\{\F_n(x)-\hat F_n(x)\right\}\,dF_0(x)= ED_n^{F_0}+O\left(\frac{\log n}{n}\right),
\end{align*}
where $D_n^{F_0}$ is defined by (\ref{def_D_n_F_0}). This is seen in the following way. Since we assume that $F_0(a)<1$, we have by by Chernoff's theorem (as in the proof of Lemma \ref{lemma:logs-exit}),
$$
\P\left\{1-\F_n(a)<\tfrac12\{1-F_0(a)\}\right\}\le e^{-nc},
$$
for a $c>0$, and hence, defining the event $A_n$ by
$$
A_n=\left\{1-\F_n(a)\ge\tfrac12\{1-F_0(a)\}\right\},
$$
we get
\begin{align*}
&\int_0^a E\left\{\F_n(x)-\hat F_n(x)\right\}\,dF_0(x)
=\int_0^a E\left\{\F_n(x)-\hat F_n(x)\right\}1_{A_n}\,dF_0(x)+O\left(e^{-nc}\right)\\
&=\int_0^a E\left\{e^{-\hat H_n(x)}-e^{-\H_n(x)}\right\}1_{A_n}\,dF_0(x)+O\left(e^{-nc}\right)\\
&=\int_0^a E\left\{1-\F_n(x)\right\}\left\{e^{-\left\{\hat H_n(x)-\H_n(x)\right\}}-1\right\}1_{A_n}\,dF_0(x)+O\left(e^{-nc}\right)\\
&=\int_0^a \left\{1-F_0(x)\right\}E\left\{e^{-\left\{\hat H_n(x)-\H_n(x)\right\}}-1\right\}1_{A_n}\,dF_0(x)+O\left(n^{-1}\right)\\
&=\int_0^a \left\{1-F_0(x)\right\}E\left\{H_0(x)+n^{-1/2}W\left(\frac{F_0(x))}{1-F_0(x)}\right)-C_n(x)\right\}\,dF_0(x)+O\left(\frac{\log n}{n}\right)\\
&=E D_n^{F_0}+O\left(\frac{\log n}{n}\right).
\end{align*}
The result now follows from Lemma \ref{lemma:df_limit}.\eop

Similarly as in Lemma \ref{lemma:remainder term}, we have the following expansion for the expectation in Corollary \ref{cor:df_limit}. The proof proceeds along the same lines as the proof of \ref{lemma:remainder term} and is therefore omitted.

\begin{lemma}
\label{lemma:remainder term2}
Assume that the conditions of Corollary \ref{cor:df_limit} are satisfied and assume in addition that $h_0$ has a bounded second derivative on $[0,a]$. Then
\begin{equation}
\label{sum_local_expectations_refined2}
n^{2/3}ED_n'
=E|C(0)|\int_0^a\left(\frac{2 h_0(t)f_0(t)}{h_0'(t)}\right)^{1/3}\,dF_0(t)+o\left(n^{-1/6}\right).
\end{equation}
\end{lemma}

The preceding results finally yield the following theorems.

\begin{theorem}
\label{th:hazard_theorem}
Let $D_n$ be defined as in Theorem \ref{th:BB_limit} and let the conditions of Theorem \ref{th:BB_limit} be satisfied.
Then:
\begin{align*}
n^{5/6}\left\{\int_0^a\left\{\H_n(x-)-\hat H_n(x)\right\}\,d\F_n(x)-ED_n\right\}\stackrel{{\cal D}}\longrightarrow N(0,\s_{H_0}^2),\,n\to\infty,
\end{align*}
where $ED_n$ and $\s_{H_0}^2$ are defined as in Theorem \ref{th:BB_limit}.
\end{theorem}

\begin{theorem}
\label{th:df_theorem}
Let the conditions of Lemma \ref{lemma:df_limit} be satisfied and let $\s_{F_0}^2$ and $C(0)$ be defined as in Lemma \ref{lemma:df_limit}. Then:
\begin{align*}
n^{5/6}\left\{\int_0^a\left\{\F_n(x-)-\hat F_n(x)\right\}\,d\F_n(x)-\int_0^a E\left\{\F_n(x-)-\hat F_n(x)\right\}\,dF_0(x)\right\}\stackrel{{\cal D}}\longrightarrow N(0,\s_{F_0}^2),\,n\to\infty,
\end{align*}
and
$$
\int_0^a E\left\{\F_n(x-)-\hat F_n(x)\right\}\,dF_0(x)\sim n^{-2/3}E|C(0)|\int_0^a\left(\frac{2h_0(t)f_0(t)}{h_0'(t)}\right)^{1/3}dF_0(t),\,n\to\infty.
$$
\end{theorem}

We only prove Theorem \ref{th:hazard_theorem}, since the proof of Theorem \ref{th:df_theorem} proceeds along similar lines.

\vspace{0.3cm}
\noindent
{\bf Proof of Theorem \ref{th:hazard_theorem}.}
Using Lemma \ref{lemma:marshall2} again, we get:
$$
\int_0^a\bigl\{\H_n(x-)-\hat H_n(x)\bigr\}\,d\F_n(x)\stackrel{{\cal D}}=\int_0^a\bigl\{H_0(x)+n^{-1/2}W_n\left(\frac{F_0(x)}{1-F_0(x)}\right)-C_n(x)\bigr\}\,d\F_n(x)+O_p\left(\frac{\log n}{n}\right),
$$
where $C_n$ is the greatest convex minorant of $V_n$ which is defined as in (\ref{V_n}) with $W_n$ replacing $W$. The process $W_n$ is distributed as standard Brownian motion on $[0,a]$ and $W_n\circ (F_0/(1-F_0))$ is given
by
$$
W_n\left(\frac{F_0(x)}{1-F_0(x)}\right)=\frac{B_n(F_0(x))}{1-F_0(x)},\,x\in[0,a],
$$
where $B_n$ is coupled to the empirical process as in Lemma \ref{embedding_lemma}.

We only have to show
\begin{equation}
\label{emp_fixed}
\int_0^a\left\{V_n(x)-C_n(x)\right\}\,d\left(\F_n-F_0\right)(x)=
o_p\left(n^{-5/6}\right),
\end{equation}
since we then have:
\begin{align*}
&\int_0^a\left\{\H_n(x-)-\hat H_n(x)\right\}\,d\F_n(x)=\int_0^a\left\{V_n(x)-C_n(x)\right\}\,d\F_n(x)+O_p\left(\frac{\log n}{n}\right)\\
&=\int_0^a\left\{V_n(x)-C_n(x)\right\}\,dF_0(x)+\int_0^a\left\{V_n(x)-C_n(x)\right\}\,d\left(\F_n-F_0\right)(x)+O_p\left(\frac{\log n}{n}\right)\\
&=\int_0^a\left\{V_n(x)-C_n(x)\right\}\,dF_0(x)+o_p\left(n^{-5/6}\right).
\end{align*}
To show that this relation holds, we follow a method which is somewhat similar to the method used in \mycite{vladrik:08} (but uses the Brownian motion representation instead of the empirical process and does not bring the derivative of the greatest convex minorant into play).

The $p$-variation of a function $f$ on the interval $I=[0,a]$
is defined by
$$
\nu_p(f;I)=\sup\left\{\sum_{i=1}^m\left|f(x_i)-f(x_{i-1})\right|^p:x_0=0<x_1<\dots<x_m=a\right\}.
$$
The $p$-variation norm of $f$ on $I$ is defined by
$$
\|f\|_{[p]}=\nu_p(f;I)^{1/p}+\sup_{x\in I}|f(x)|.
$$
We have, by Theorem II.3.27 in \mycite{dudley:99}, for $p,q>0$ and $1/p+1/q>1$:
\begin{equation}
\label{p_norm_integral}
\left|\int_{[0,a]}\left\{V_n(x)-C_n(x)\right\}\,d\left(\F_n-F_0\right)(x)\right|
\le c\,\|V_n-C_n\|_{[p]}\|\F_n-F_0\|_{[q]},
\end{equation}
for a constant $c>0$. Moreover, by Theorems I.6.1 and I.6.2 in \mycite{dudley:99}, and Theorem 3.2 in \mycite{qian:98}, we have:
\begin{equation}
\label{p-norm_emp}
\|\F_n-F_0\|_{[q]}=\left\{\begin{array}{lll}
O_p\left(n^{(1-q)/q}\right)\,&,\,q\in[1,2),\\
O_p\left(n^{-1/2}\sqrt{L(Ln))}\right)\,&,\,q=2,\\
O_p\left(n^{-1/2}\right)\,&,\,q>2,
\end{array}
\right.
\end{equation}
where $Ln=1\vee\log n$.

Let $\t_1,\dots,\t_m$ be the points of jump of the derivative $c_n$ of $C_n$ on $[0,a]$, and let $\t_0=0,\,\t_{m+1}=a$.
The function $C_n$ is linear on the intervals $[\t_i,\t_{i+1}]$, and $V_n$ behaves on such an interval as an excursion above its greatest convex minorant $C_n$, with the same values as $V_n$ at the endpoints of the interval. Hence we have, for $p>2$, by Lemma 4 of \mycite{huang_dudley:01},
\begin{align*}
&\nu_p(V_n-C_n;[0,a])\le 2^{p-1}\sum_{k=1}^{m+1}\nu_p(V_n-C_n;[\t_{i-1},\t_i])
=2^{p-1}\sum_{i=1}^{m+1}\nu_p(V_n;[\t_{i-1},\t_i]),
\end{align*}
where, using the fact that the linear part drops out in taking the comparison with the greatest convex minorant,
\begin{align*}
\tilde V_n(x)&=n^{-1/2}\left\{W\left(\frac{F_0(x)}{1-F_0(x)}-\frac{F_0(\t_{i-1})}{1-F_0(\t_{i-1})}\right)
-\frac{x-\t_{i-1}}{\t_i-\t_{i-1}}W\left(\frac{F_0(\t_i)}{1-F_0(\t_i)}-\frac{F_0(\t_{i-1})}{1-F_0(\t_{i-1})}\right)\right\}\\
&\qquad\qquad\qquad\qquad+H_0(x)-H_0(\t_{i-1})-\frac{x-\t_{i-1}}{\t_i-\t_{i-1}}\left\{H_0(\t_i)-H_0(\t_{i-1})\right\},\,x\in[\t_{i-1},\t_i].
\end{align*}
By part (ii) of Lemma \ref{localization_lemma} we have:
$$
E\max_i(\t_i-\t_{i-1})=O\left(n^{-1/3}\log n\right).
$$
Let $u_i$ be the midpoint of the interval $[\t_{i-1},\t_i]$ and let $f_{H_0}$ by defined by
$$
f_{H_0}(x)=H_0(x)-H_0(\t_{i-1})-\frac{x-\t_{i-1}}{\t_i-\t_{i-1}}\left\{H_0(\t_i)-H_0(\t_{i-1})\right\},\,x\in\left[\t_{i-1},\t_i\right].
$$
Then
\begin{align*}
f_{H_0}(x)=-\tfrac12h_0'(u_i)\{x-\t_{i-1}\}\{\t_i-x\}\{1+o_p(1)\},
\end{align*}
where $x\mapsto \{x-\t_{i-1}\}\{\t_i-x\}$ is increasing on $[\t_{i-1},u_i]$ and decreasing on $[u_i,\t_i]$, and
$$
\nu_p\left(f_{H_0};[\t_{i-1},\t_i])\right)\sim 2^{1-p}h_0'(u_i)^p\left\{u_i-\t_{i-1})\right\}^p\left\{\t_i-u_i)\right\}^p,
$$
(see, e.g., (3.4) of \mycite{huang_dudley:01}). Hence, for any $p>2$,
\begin{align*}
&\sum_{i=1}^{m+1}\nu_p(f_{H_0};[\t_{i-1},\t_i])\\
&\sim2^{1-p}\sum_{i=1}^{m+1}h_0'(u_i)^p\left\{u_i-\t_{i-1}\right\}^p\left\{\t_i-u_i\right\}^p
=2^{1-3p}\sum_{i=1}^{m+1}h_0'(u_i)^{p}\left\{u_i-\t_{i-1}\right\}^{2p}\\
&\le2^{-3p}\max_i\left\{u_i-\t_{i-1}\right\}^{2p-1}\sum_{i=1}^{m+1}h_0'(u_i)^{p}\left\{\t_i-\t_{i-1}\right\}\\
&\sim2^{-5p+1}\max_i\left\{\t_i-\t_{i-1}\right\}^{2p-1}\int_0^a h_0'(u)^p\,du
=O_p\left(n^{-(2p-1)/3}(\log n)^{(2p-1)/2}\right).
\end{align*}
Note that the $O_p$-term becomes $O_p\left(n^{-1}(\log n)^{3/2}\right)$ for $p=2$.

For the Brownian part
$$
B_{nk}(x)\stackrel{\mbox{\small def}}=n^{-1/2}\left\{W\left(\frac{F_0(x)}{1-F_0(x)}-\frac{F_0(\t_{i-1})}{1-F_0(\t_{i-1})}\right)
-\frac{x-\t_{i-1}}{\t_i-\t_{i-1}}W\left(\frac{F_0(\t_i)}{1-F_0(\t_i)}-\frac{F_0(\t_{i-1})}{1-F_0(\t_{i-1})}\right)\right\},
$$
we find, for $p>2$,
$$
\sum_{i=1}^{m+1}\nu_p(B_{nk};[\t_{i-1},\t_i])=O_p\left(n^{-p/2}\right),
$$
by the fact that almost all Brownian motion paths are H\"older continuous of any order $<1/2$.

So we find:
\begin{equation}
\label{p_norm_integrand}
\|V_n-C_n\|_{[p]}=O_p\left(n^{-1/2}(\log n)^{(2p-1)/(2p)}\right),
\end{equation}
for any $p>2$. Thus (\ref{p_norm_integral}), (\ref{p-norm_emp}) and (\ref{p_norm_integrand}) imply
$$
\int_{[0,a]}\left\{V_n(x)-C_n(x)\right\}\,d\left(\F_n-F_0\right)(x)=O_p\left(n^{-1+\e}\right),
$$
for arbitrarily small $\e>0$.\eop

\vspace{0.3cm}
We end this section with a result for the situation that the hazard is nondecreasing, but not strictly nondecreasing.

\begin{theorem}
\label{th:df_theorem2}
Let $\hat F_n$ and $\F_n$ be defined as in Theorem \ref{th:df_theorem} and let (again)
$$
U_n=\int_{[0,a]}\bigl\{\F_n(x-)-\hat F_n(x)\bigr\}\,d\F_n(x).
$$
Let $U$ be given by
$$
U=\int_0^a\bigl\{1-F_0(x)\bigr\}\left\{W\left(\frac{F_0(x)}{1-F_0(x)}\right)-C(x)\right\}\,dF_0(x),
$$
where $W$ is standard Brownian motion on $[0,\infty)$ and $C$ is the greatest convex minorant of
\begin{equation}
\label{BM_for_constant_haz}
x\mapsto W\left(\frac{F_0(x)}{1-F_0(x)}\right),\,x\in[0,a].
\end{equation}
Suppose that the underlying hazard $h_0$ is constant on $[0,a]$. Then:
\begin{align*}
n^{1/2}U_n\stackrel{{\cal D}}\longrightarrow U,\,n\to\infty.
\end{align*}
\end{theorem}

\noindent
{\bf Proof.} The proof follows lines that are familiar by now. We first consider
$$
U_n'=\int_{[0,a]}\bigl\{\F_n(x-)-\hat F_n(x)\bigr\}\,dF_0(x).
$$
By (\ref{BM_replacement}) we can replace $\F_n-\hat F_n$ by:
$$
\exp\left\{-C_n(x)\right\}\left\{1-\exp\left\{-H_0(x)-n^{-1/2}W\left(\frac{F_0(x))}{1-F_0(x)}\right)-C_n(x)\right\}\right\},
$$
where $C_n$ is the greatest convex minorant of the process
$$
x\mapsto H_0(x)+n^{-1/2}W\left(\frac{F_0(x))}{1-F_0(x)}\right),\,x\in[0,a],
$$
with a remainder term of order $O_p((\log n)/n)$. Using the delta method as in the proof of Corollary \ref{cor:df_limit}, we can replace this (apart from a remainder term of order $O_p(n^{-1})$) by:
$$
n^{-1/2}\{1-F_0(x)\}\left\{W\left(\frac{F_0(x))}{1-F_0(x)}\right)-C(x)\right\},\,x\in[0,a].
$$
where $C$ is the greatest convex minorant of the process (\ref{BM_for_constant_haz}), and where we use that $H_0$ is linear on $[0,a]$. The statement for $U_n$ now follows by an application of \mycite{dudley:99}, as in the proof of Theorem \ref{th:hazard_theorem}.
\eop

\begin{remark}
{\rm Note that the rate of convergence drops from $n^{5/6}$ to $n^{1/2}$ in Theorem \ref{th:df_theorem2}, and that the limiting distribution is not normal. We get a limit behavior that can be analyzed using the methods
of \mycite{piet:83}, where the concave majorant of Brownian motion without drift is characterized via a Poisson process of jump locations and Brownian excursions.}
\end{remark}

\section*{Appendix}
\noindent
{\bf Proof of Lemma \ref{lem:Vbelow0}.}
Let $u>0$. Then, for $x\ge u$:
$$
V(x)=W(x)+(x-u)^2+2u(x-u)+u^2\ge W(x)+(x-u)^2+u^2=W(u)+u^2+W(x)-W(u)+(x-u)^2.
$$
Hence,
\begin{align*}
&\P\left(\min_{x\ge u}V(x)\le0\right)\le \P\left(\min_{x\ge u}W(u)+u^2+W(x)-W(u)+(x-u)^2\le0\right)=\\
&\P\left(W(u)+u^2+\min_{x\ge u}W(x)-W(u)+(x-u)^2\le0\right)\le \\& \P\left(W(u)\le -\tfrac12u^2\right)+\P\left(\min_{x\ge u}W(x)-W(u)+(x-u)^2\le-\tfrac12 u^2\right)
\end{align*}
The process
$$
x\mapsto W(x)-W(u)+(x-u)^2,\,x\ge u
$$
behaves in the same way as the process $t\mapsto V(t),\,t\ge0$, but starts in $x$ instead of $0$.
By Corollary 2.1 in \mycite{piet_nico:10} we have that for all  $z>0$,
\begin{equation}
\label{Airy_asymp}
\P\left\{\min_{t\in\R}V(t)\le -z\right\}\sim 2\cdot3^{-1/2}\exp\left\{-8z^{3/2}/\sqrt{27}\right\},\,z\to\infty,
\end{equation}
implying that there exist positive constants $c_1$ and $c_2$ such that for all $u\ge 0$
$$
\P\left(\min_{x\ge u}W(x)-W(u)+(x-u)^2\le-\tfrac12 u^2\right)\le c_1 \exp\left\{-c_2u^{3}\right\}.
$$
We also have for all $u>0$
$$
\P\left\{W(u)<-\tfrac12u^2\right\}=\P\left\{W(u)/\sqrt{u}<-\tfrac12u^{3/2}\right\}\le  \frac{\exp\left\{-\tfrac1{8}u^3\right\}}{u^{3/2}\sqrt{\pi/2}}
$$
implying that there exist positive constants $c_3$ and $c_4$ such that for all $u\ge 0$
$$
\P\left\{W(u)<-\tfrac12u^2\right\}\le c_3 \exp\left\{-c_4u^{3}\right\}.
$$
Combining these upper bounds with the fact that the process $V$ running to the left from zero behaves in the same way as the process $V$ running to the right from zero, the statement of the lemma follows.
\eop

\vspace{0.3cm}
\noindent
{\bf Proof of Lemma \ref{localization_lemma}.} (i). The interval $I_{n,k}$ is bounded on the left by the interval $\tilde{J}_{n,k}$ and on the right by the interval $\bar{J}_{n,k+1}$. The intervals $\tilde{J}_{n,k}$ and $\bar{J}_{n,k+1}$ both have length of order $n^{-1/3}\sqrt{\log n}$. If the greatest convex minorant $C_n$ of $V_n$ on $[0,a]$ has changes of slope in the intervals $\tilde{J}_{n,k}$ and $\bar{J}_{n,k+1}$, the greatest convex minorant of $V_n$ on $[0,a]$, restricted to the interval $I_{n,k}$, coincides with the greatest convex minorant $C_{nk}$ of $V_n$ on $L_{n,k}$, restricted to the interval $I_{n,k}$. So we have to find bounds for the probability that the greatest convex minorant of $V_n$ on $[0,a]$ has no changes of slope in $\tilde{J}_{nk}$ or $\bar{J}_{n,k+1}$. To do this, we  follow the method used in \mycite{grwe:92}, p.\ 96.

Let $a_{nk}$ and $b_{nk}$ be the left and right endpoints of $\bar{J}_{n,k+1}$, respectively, and let $u_{nk}$ be its midpoint. If
\begin{equation}
\label{assumption00}
c_n(a_{nk})<h_0(u_{nk})<c_n(b_{nk}),
\end{equation}
where $c_n$ is the left-continuous slope of $C_n$, then $C_n$ has a change of slope in the interval $\bar{J}_{n,k+1}$. Note that for $x\ge b_{nk}$, using the assumed smoothness of $H_0$, and $\inf_{[0,a]}h_0^\prime(x)=2\kappa>0$,
\begin{equation}\label{eq:Vnlower} V_n(x)-V_n(u_{nk})\ge n^{-1/2}\left\{W\left(\frac{F_0(x)}{1-F_0(x)}\right)-W\left(\frac{F_0(u_{nk})}{1-F_0(u_{nk})}\right)\right\}+h_0(u_{nk})(x-u_{nk})+\kappa (x-u_{nk})^2.\end{equation}

Now consider the event that
\begin{equation}
\label{assumption1}
c_n(b_{nk})\le h_0(u_{nk}),
\end{equation}
and let $\t_{nk}$ be the first point of jump of $c_n$ to the right of $b_{nk}$. Then
$$
c_n(x)\le h_0(u_{nk}),\,x<\t_{nk},
$$
and hence
$$V_n(\t_{nk})-V_n(x)
\le C_n(\t_{nk})-C_n(x)=\int_x^{\t_{nk}}c_n(y)\,dy\le h_0(u_{nk})(\t_{nk}-x),\,x<\t_{nk}.
$$
Using (\ref{eq:Vnlower}) and stationarity of Brownian Motion, this means that the probability of (\ref{assumption1}) is bounded above by
\begin{align}
&\P\left\{V_n(\t_{nk})-V_n(u_{nk})\le h_0(u_{nk})(\t_{nk}-u_{nk})\right\}\le\P\left\{\exists x\ge b_{nk}\,:\,V_n(x)-V_n(u_{nk})\le h_0(u_{nk})(x-u_{nk})\right\}\nonumber\\
&\qquad \le \P\left\{\exists x\ge b_{nk}\,:\,    n^{-1/2}\left\{W\left(\frac{F_0(x)}{1-F_0(x)}\right)-W\left(\frac{F_0(u_{nk})}{1-F_0(u_{nk})}\right)\right\}\le -\kappa (x-u_{nk})^2 \right\}\nonumber\\
&\qquad=\P\left\{\exists x\ge b_{nk}\,:\,    n^{-1/2}\left\{W\left(\frac{F_0(x)}{1-F_0(x)}-\frac{F_0(u_{nk})}{1-F_0(u_{nk})}\right)\right\}\le -\kappa (x-u_{nk})^2 \right\}.\label{eq:probVn}
\end{align}
We will see that this probability will become exponentially small. To this end, define the following covering of $[b_{nk},a]$
$$
K_{nkj}\stackrel{\mbox{\small def}}=\left[t_{nkj},t_{nk,j+1}\right]\stackrel{\mbox{\small def}}=\left[b_{nk}+jn^{-1/3},b_{nk}+(j+1)n^{-1/3}\right]=[b_{nk},a]
$$
for $0\le j \le \lfloor n^{1/3}(a-b_{nk})\rfloor$ (where the right end point of the last interval is taken to be $a$). Then the probability in (\ref{eq:probVn}) can be bounded above by
\begin{equation}
\label{eq:sumbound}
\sum_{j=0}^{\lfloor n^{1/3}(a-b_{nk})\rfloor}\P\left\{\exists x\in K_{nkj}\,:\,n^{-1/2}\left\{W\left(\frac{F_0(x)}{1-F_0(x)}-\frac{F_0(u_{nk})}{1-F_0(u_{nk})}\right)\right\}\le -\kappa (x-u_{nk})^2 \right\}.
\end{equation}
Denoting the probabilities in this sum by $p_{nkj}$, we get
\begin{align*}
&p_{nkj}\le\P\left\{\sup_{ x\in K_{nkj}}W\left(\frac{F_0(x)}{1-F_0(x)}-\frac{F_0(u_{nk})}{1-F_0(u_{nk})}\right)
\ge \kappa \sqrt{n}(t_{nkj}-u_{nk})^2 \right\}\\
&\qquad \le \P\left\{\sup_{ 0\le z\le F_0(t_{nk,j+1})/(1-F_0(t_{nk,j+1}))-F_0(u_{nk})/(1-F_0(u_{nk}))}W(z)
\ge \kappa \sqrt{n}(t_{nkj}-u_{nk})^2 \right\}.
\end{align*}
Since $t_{nk,j+1}\in[b_{nk},a]$ for all $j$ under consideration,
$$
0\le \frac{F_0(t_{nk,j+1})}{1-F_0(t_{nk,j+1})}-\frac{F_0(u_{nk})}{1-F_0(u_{nk})}\le \frac{(F_0(t_{nk,j+1})-F_0(u_{nk}))}{(1-F_0(a))^2}\le \lambda(t_{nk,j+1}-u_{nk})
$$
for some $0<\lambda<\infty$, we obtain, for a standard normal random variable $Z$
\begin{align*}
&p_{nkj}\le \P\left\{\sup_{ 0\le z\le \lambda(t_{nk,j+1}-u_{nk})}W(z)
\ge \kappa \sqrt{n}(t_{nkj}-u_{nk})^2 \right\}=\P\left\{|Z|\ge\frac{\kappa \sqrt{n}(t_{nkj}-u_{nk})^2}{\sqrt{ \lambda(t_{nk,j+1}-u_{nk})}}\right\}\\
&\qquad\le\P\left\{|Z|\ge \tilde{\kappa}\sqrt{n}(t_{nkj}-u_{nk})^{3/2}\right\}\le \frac12\exp\left\{-\tfrac12n\tilde{\kappa}^2(t_{nkj}-u_{nk})^{3}\right\}.
\end{align*}
Using that $t_{nkj}-u_{nk}=b_{nk}-u_{nk}+jn^{-1/3}$, and $b_{nk}-u_{nk}\sim \tfrac12n^{-1/3}\sqrt{\log n}$, we get
$$
p_{nkj}\le \exp\left\{-\tfrac12\tilde{\kappa}  \left((\log n)^{3/2}+j^3\right)\right\}\Rightarrow \sum_{j=0}^{\lfloor n^{1/3}(a-b_{nk})\rfloor}p_{nkj}\le  \rho\exp\left\{-\rho^\prime(\log n)^{3/2}\right\}
$$
for some $\rho,\rho^\prime>0$. Combining this with (\ref{eq:probVn}) and (\ref{eq:sumbound}), this bounds the probability of (\ref{assumption1}) from above. Since a similar bound holds for the probability of the event $c_n(a_{nk})\ge h_0(u_{nk})$, the probability that (\ref{assumption00}) does not hold for a specific $k$, is bounded by a bound of the same structure. Moreover, since this upper bound does not depend on $k$ and  $m_n\sim an^{1/3}/\log n$, the probability that there exists a $1\le k\le m_n$ for which   (\ref{assumption00})  does not hold satisfies the same bound (with slight change in $\rho$ and $\rho^\prime$), this proves (i). Part (ii) is an immediate consequence of (i).
\eop

\end{document}

%% file: commands3.tex
\newcommand{\mycite}[1]{{\small \sc \citeNP{#1}}}

\def\R{\mathbb R}

\def\F{{\mathbb F}}

\def\H{{\mathbb H}}

\def\P{{\mathbb P}}

\def\labda1{\lambda_1}
\def\labda2{\lambda_2}

\def\e{\varepsilon}
\def\f{\phi}
\def\t{\tau}

\def\s{\sigma}

\def\comment#1{\relax}

\def\=in{\mathop{\rm =}}

\def\eop{\hfill\mbox{$\Box$}\newline}

\newtheorem{theorem}{Theorem}[section]
\newtheorem{corollary}{Corollary}[section]
\newtheorem{lemma}{Lemma}[section]

\newtheorem{remark}{Remark}[section]